%% file: 3rank.tex
\documentclass{article}
\pdfoutput=1
\usepackage{amssymb,amsfonts,amsmath,amsthm}
\usepackage{epsfig}
\usepackage{longtable}
\usepackage[section]{algorithm}
\usepackage[noend]{algorithmic}

\include{custommacros}

\title{Computing quadratic function fields with high $3$-rank via cubic field tabulation}

\author{Pieter Rozenhart \\ Department of Mathematics and Statistics, University of
  Calgary,\\
  2500 University Drive NW, Calgary, Alberta, Canada, T2N 1N4  \\ \texttt{pmrozenh@alumni.uwaterloo.ca} \and Michael Jacobson Jr.\ \\   Department of Computer Science, University of
  Calgary,\\
  2500 University Drive NW, Calgary, Alberta, Canada, T2N 1N4\\
   \texttt{jacobs@ucalgary.ca} 
\and Renate Scheidler \\Department of Mathematics and Statistics, University of
  Calgary,\\
  2500 University Drive NW, Calgary, Alberta, Canada, T2N 1N4\\
    \texttt{rscheidl@ucalgary.ca}}

\begin{document}
\maketitle

\begin{abstract}
This paper presents an algorithm for generating all imaginary and unusual discriminants up to a fixed degree bound that define a quadratic function field of positive 3-rank. Our method makes use of function field adaptations of a method due to Belabas for finding quadratic number fields of high 3-rank and of a refined function field version of a theorem due to Hasse. We provide numerical data for discriminant degree up to 11 over the finite fields $\finfldq{5}, \finfldq{7}, \finfldq{11}$ and $\finfldq{13}$. A special feature of our technique is that it produces quadratic function fields of minimal genus for any given 3-rank. Taking advantage of certain $\ratsfuncq{q}$-automorphisms in conjunction with Horner's rule for evaluating polynomials significantly speeds up our algorithm in the imaginary case; this improvement is unique to function fields and does not apply to number field tabulation. These automorphisms also account for certain divisibility properties in the number of fields found with positive $3$-rank.  Our numerical data mostly agrees with the predicted heuristics of Friedman-Washington and partial results on the distribution of such values due to Ellenberg-Venkatesh-Westerland for quadratic function fields over the finite field $\finfldq{q}$ where $q \equiv -1 \pmod{3}$.  The corresponding data for  $q \equiv 1 \pmod{3}$ does not agree closely with the previously mentioned heuristics and results, but does agree more closely with some recent number field conjectures of Malle and some work in progress on proving such conjectures for function fields due to Garton. \end{abstract}

\section{Introduction}\label{sec:intro}
Let $D$ be a square-free non-constant polynomial in \polyringq{q} and $\mbox{Cl}(D)$ the ideal class group of the quadratic function field $\ratsfuncqextn{q}{\sqrt{D}}$.  For any prime $\ell$, the number $r_{\ell}(D)$, which denotes the number of cyclic factors in the $\ell$-Sylow subgroup of $\mbox{Cl}(D)$, is called the \emph{$\ell$-rank} of the group $\mbox{Cl}(D)$. In short, we say that the quadratic function field $\ratsfuncqextn{q}{\sqrt{D}}$ has $\ell$-rank $r$ if $\mbox{Cl}(D)$ has $\ell$-rank equal to $r$.  

In this paper, we develop an algorithm for finding all quadratic discriminants $D$ (with $-3D$ imaginary or unusual) of bounded degree for which $\mbox{Cl}(D)$ has positive $3$-rank.
Our algorithm is based on our previous work (\cite{ANTSPie}, \cite{Pieter2} and \cite{Pieter3}) for tabulating cubic function fields of bounded degree, and
is inspired by Belabas' algorithm \cite{Bela2} for tabulating quadratic number fields
of bounded discriminant with three rank greater than zero.  Like Belabas' algorithm, our algorithm makes use of an extended function field version of a theorem of Hasse \cite{Hasse1} relating the $3$-rank of a quadratic function field of discriminant $D$ to the number of triples of conjugate cubic function fields with discriminant $D$.  
This theorem is also used in \cite{cuffqi}, although it is used in the ``reverse" direction, in the sense that information on the ideal class group of quadratic function fields is used to generate cubic function fields.  Our approach is the opposite:  we tabulate cubic function fields and use this data to generate information on $3$-ranks of quadratic fields.  We also make use of the $\ratsfuncq{q}$-automorphisms that send $t$ to $t+\alpha$ with $\alpha \in \finfldq{q}^*$ to speed up our algorithm; this new improvement is unique to the function field setting, and results in a speedup by a factor of approximately $q$ in most cases over a basic version of our algorithm.  These $\ratsfuncq{q}$-automorphisms also explain why the number of fields with a given positive 3-rank has certain divisibility properties.

While the basic ideas of this paper stem from Chapter 7 of the first author's doctoral dissertation \cite{Pieter2}, written under the supervision of the last two authors, a large portion of the material herein is new and represents a significant improvement to \cite{Pieter2}.
Our contributions are as follows.  First, our basic algorithm (Algorithm \ref{tabimag}) functions in much the same way as that in \cite{Pieter3}:  we loop over the coefficients of a reduced binary cubic form except instead of outputting a binary cubic form of discriminant $D$, we increment a counter for that discriminant value and store both $D$ and its counter value in a hash table.  The hash table is then output once the for loops are exited.  The improved algorithm for the imaginary case (Algorithm \ref{tabimagHorner}) shortens the loops on two of the coefficients.  This results in some forms being missed, and hence some of the discriminant counters are short of their actual values.  In order to recover the correct values for the discriminant counters, we check which discriminants had forms of a certain type, and augment those discriminant counters as necessary.  Complete details appear in Sections \ref{sec:alg} and \ref{subsec:newalg}.  The improvement to the algorithm is unique to the function field setting, and results in asymptotic improvements by a factor of $q$ in most cases.  These improvements allowed us to push our computations beyond those in \cite{Pieter2}.

Our $3$-rank algorithm, like Belabas' \cite{Bela2}, is exhaustive in the sense that all quadratic discriminants of any fixed degree that define quadratic function fields with positive $3$-rank are produced by the method.  As a consequence, the resulting algorithm also produces minimum discriminant sizes for any given positive $3$-rank value.  This is in contrast to \cite{BauJac} and \cite{RenateWIN}, which present methods for finding quadratic function fields with high 3-rank whose discriminants are small, but not necessarily minimal.  Our goal is to obtain these minimal discriminant sizes for a given positive $3$-rank value $r$, rather than record $3$-rank values. 

Our second contribution is that since our method is exhaustive, we can generate data on the distribution of $3$-rank values up to a fixed bound $B$ on $\deg(D)$.  Cohen and Lenstra \cite{CohLenHeur1,CohLenHeur2} gave heuristics on the behavior of class groups of quadratic number fields.  For example, they provided heuristic estimates for the probability that the $\ell$-rank of a given class group is equal to $r$ for a given prime $\ell$ and non-negative integer $r$. None of these heuristics are proved, but there is a large amount of numerical evidence supporting their validity, as seen, for example, in \cite{JacRamaWill,teRieleWill1}.   These heuristics imply that the ideal class group of a quadratic number field is expected to have low $\ell$-rank for any prime $\ell$.  Consequently, there is a large body of literature devoted to the construction of families of quadratic number fields of large $\ell$-rank, with $3$-ranks of particular interest.  

The function field analogue of the Cohen-Lenstra heuristics, the Friedman-Washington heuristics \cite{FriedWash}, attempt to explain statistical observations about divisor class groups of quadratic function fields.   Some progress has been made in trying to prove these heuristics, most notably by Achter \cite{Achter1,Achter2}, Ellenberg, Venkatesh and Westerland \cite{EllenVenkWest,EllenVenk3} and Garton \cite{Garton}.  These results are somewhat weaker than the original heuristics, as they rely on the size $q$ of the underlying finite field tending to infinity in addition to the genus of the function field.  Previous results attempting to numerically verify the Friedman-Washington heuristics include computations of class groups for small genus over small base fields by Feng and Sun \cite{FengSun1} and computations of class groups of real quadratic function fields of genus $1$ over large base fields by Friesen \cite{Fries1}.  To the knowledge of the authors, the computational data contained herein is the most extensive since the work of Feng and Sun \cite{FengSun1} and Friesen \cite{Fries1}.

We were able to generate examples of minimal genus and 3-rank values as large as four for quadratic function fields over $\finfldq{q}$ for $q=5,7,11,13$.  As expected, we did not find fields with higher 3-rank than any known examples, but we did find numerous examples of fields with 3-rank as high as four and smaller genus than any others known.  In addition, the data we generated yields evidence for the validity of the Friedman-Washington heuristic for $q=5,11$.  Due to the presence of cube roots of unity, the data for $q=7,13$ does not agree closely with Friedman-Washington/Ellenberg et al., but instead with some recent conjectures of Malle \cite{Malle,Malle2} for number fields, and a new distribution result for function fields of Garton \cite{Garton,Garton2}.  Our data suggests that all the doubly-asymptotic results of Ellenberg et al.\ and Garton, where $q \rightarrow \infty$ as well as $g \rightarrow \infty$, may also hold for fixed $q$.  
%

This paper is organized as follows.  After a brief review of some preliminaries from the theory of algebraic function fields and cubic function field tabulation in Section \ref{sec:prel}, we proceed with a short discussion of the basic algorithm in Section \ref{sec:alg}, and the improved algorithm for imaginary discriminants in Section \ref{subsec:newalg}.  The algorithm's complexity is analyzed in Section \ref{sec:algcomplexity}.  We discuss the Friedman-Washington heuristics and related results in Section \ref{sec:3rankbg}.  The $3$-rank data generated is presented in Section  \ref{sec:3ranknnumericalresultsimag}, and a comparison to the Friedman-Washington heuristics, the Achter and Ellenberg et al.\ distribution, Malle's conjectured formula and Garton's distribution are presented in Section \ref{ssec:AutoCounts}.  Finally, we make some concluding remarks and suggest open problems in Section \ref{sec:concl}.

\section{Preliminaries}\label{sec:prel}

Let \finfldq{q}\ be a finite field of characteristic at least~$5$, and
set $\finfldq{q}^* = \finfldq{q} \backslash \{0\}$. Denote by
\polyringq{q}\ and \ratsfuncq{q}\ the ring of polynomials and the
field of rational functions in the variable $t$ over \finfldq{q},
respectively.  For any non-zero $H \in \polyringq{q}$ of degree $n
= \deg(H)$, we let $|H| = q^n = q^{\deg(H)}$, and denote by
$\sgn(H)$ the leading coefficient of $H$. For $H = 0$, we set $|H|
= 0$. This absolute value extends in the obvious way to
\ratsfuncq{q}. Note that in contrast to the absolute value on the
rational numbers \rats, the absolute value on \ratsfuncq{q} is non-Archimedean.



A \emph{binary quadratic form over \polyringq{q}} is a homogeneous quadratic polynomial in two variables with coefficients in
\polyringq{q}.  We denote the binary quadratic form $H(x,y) = Px^2 + Qxy + Ry^2$ by  $H = (P,Q,R)$.
The \emph{discriminant} of $H$ is the polynomial
$D(H) = Q^2 - 4PR \in \polyringq{q}$.  A polynomial $F$ in \polyringq{q}\ is said to be
\emph{imaginary} if $\deg(F)$ is odd, \emph{unusual} if
$\deg(F)$ is even and $\sgn(F)$ is a non-square in
$\finfldq{q}^*$, and \emph{real} if $\deg(F)$ is even and
$\sgn(F)$ is a square in $\finfldq{q}^*$.  Correspondingly, a binary quadratic form is said to be \emph{imaginary, unusual or real} according to whether its discriminant is imaginary, unusual or real.

A \emph{binary cubic form over \polyringq{q}} is a homogeneous
cubic polynomial in two variables with coefficients in
\polyringq{q}.  We denote the binary cubic form $f(x,y) = ax^3 + bx^2y + cxy^2 + dy^3$ by $f =
(a,b,c,d)$.  The \emph{discriminant} of $f = (a,b,c,d)$
is the polynomial
\[ D(f) = 18abcd + b^2c^2 - 4ac^3 - 4b^3d -
27a^2d^2 \in \polyringq{q} \enspace . \]
We assume throughout that binary cubic forms $f = (a,b,c,d)$ are \emph{primitive}, i.e. $\gcd(a,b,c,d) =
1$, and irreducible.  

Let $F$ be a binary quadratic or cubic form over \polyringq{q}. If $M = \left(\begin{smallmatrix} \alpha & \beta \\
    \gamma & \delta \end{smallmatrix} \right)$
is a $2 \times 2$ matrix with entries in \polyringq{q}, then the
\emph{action} of $M$ on $F$ is defined by $(F \circ M)(x,y) = F(\alpha x
+ \beta y, \gamma x + \delta y)$.  We obtain an equivalence relation from this action by restricting
to matrices $M \in \genlinfunc$, the group of $2 \times 2$
matrices over $\polyringq{q}$ whose determinant lies in
$\finfldq{q}^*$. That is, two binary quadratic or cubic forms $F$
and $G$ over \polyringq{q} are said to be \emph{equivalent} if
\[F(\alpha x + \beta y, \gamma x + \delta y) = G(x,y)\]
for some $\alpha, \beta, \gamma,
\delta \in \polyringq{q}$ with $\alpha \delta - \beta \gamma \in
\finfldq{q}^*$. Up to some even power of $\det(M)$, equivalent binary forms have the
same discriminant. Furthermore, the action of the group
$\genlinfunc$ on binary forms over \polyringq{q} preserves
irreducibility and primitivity over~\ratsfuncq{q} .

As in the case of integral binary cubic forms, any binary cubic
form $f = (a,b,c,d)$ over \polyringq{q} is closely associated with
its {\em Hessian}
\renewcommand{\arraystretch}{2}
\[ H_f(x,y) = -\frac{1}{4}\left|\begin{array}{cc}
\displaystyle\frac{\partial^2f}{\partial x \partial x} &
\displaystyle\frac{\partial^2f}{\partial x \partial y} \\
\displaystyle\frac{\partial^2f}{\partial y \partial x} &
\displaystyle\frac{\partial^2f}{\partial y \partial y}\end{array}
\right| = (P,Q,R) \enspace ,\]
where $P = b^2 - 3ac$, $Q = bc - 9ad$, and $R = c^2 - 3bd$. Note
that $H_f$ is a binary quadratic form over \polyringq{q}. The
Hessian has a number of useful properties, which are easily
verified by direct computation:

\renewcommand{\arraystretch}{1}

\[H_{f \circ M} = (\det   M)^2(H_f \circ M)\mbox{ for any }M \in \genlinfunc, \mbox{ and}\]
\[D(H_f) = -3\, D(f).\]

We now briefly summarize the reduction theory for binary quadratic and cubic forms.  Fix a primitive root $h$ of $\finfldq{q}^*$.  As in Artin \cite{Artin1}, we only consider quadratic discriminants $D$ endowed with the normalization $\sgn(D) = 1$ or $\sgn(D) = h$, where $1$ or $h$ is chosen depending on whether or not $\sgn(D)$ is a square in $\finfldq{q}^*$.   We can impose this restriction since the discriminant of a function field is only unique up to square factors in $\finfldq{q}^*$.  Define the set $S = \{ h^i :  0 \leq i \leq (q-3)/2 \}$.  Then $a \in S$ if and only if $-a \notin S$.  In particular, note that $S$ is non-empty, since $1 \in S$.

\begin{definition}[(Rozenhart \cite{Pieter2})]
\begin{enumerate}
\item Let $H = (P,Q,R)$ be an imaginary or unusual binary quadratic form of discriminant $D$. Then 
$H$ is \emph{reduced} if
\begin{enumerate}
\item $|Q| < |P|$, and either $Q = 0$ or $\sgn(Q) \in S$;
\item Either $|P| < |R|$ and $\sgn(P) \in \{1, h\}$, or $|P| = |R|$ and $\sgn(P) = 1$;
\item When $|P| = |R|$ and $\sgn(P)=1$, then $f$ is lexicographically smallest among the $q+1$ binary quadratic forms 
in its equivalence class satisfying conditions (a) and (b) above.
\end{enumerate}
\item Let $f = (a,b,c,d)$ be a binary cubic form with imaginary or unusual Hessian $H_f = (P, Q, R)$ of discriminant $-3D$. Then $f$ is \emph{reduced} if
\begin{enumerate}
\item $\sgn(a) \in S$, and if $Q = 0$, then $\sgn(d) \in S$.
\item $H_f$ is reduced, and in addition, if $|P| = |R|$, then $f$ is 
lexicographically smallest among all binary cubic forms in its equivalence class with 
Hessian $H_f$.
\end{enumerate}
\end{enumerate}
\label{reduceddefn}
\end{definition}

Note that there are rare occasions where two equivalent binary cubic forms can have 
the same reduced unusual Hessian $H_f = (P,Q,R)$; see Theorem 4.21 of \cite{Pieter2}. 
The proof of the following theorem can be found in \cite{Pieter2}.

\begin{theorem} \mbox{}
\begin{enumerate}
\item Every equivalence class of imaginary or unusual binary quadratic forms contains a 
unique reduced representative, and there are only finitely many reduced 
imaginary or unusual binary quadratic forms of any given discriminant.
\item Every equivalence class of binary cubic forms with imaginary or unusual Hessian 
contains a unique reduced representative, and there are only finitely 
many reduced binary cubic forms of any given discriminant with imaginary or unusual 
Hessian.
\end{enumerate}
\label{uniquenessofredforms}
\end{theorem}

Recall from Theorem 5.3 of \cite{Pieter3} that if $f=(a,b,c,d)$ is a reduced binary cubic form of discriminant $D$, then $\deg(a), \deg(b) \leq \deg(D)/4$, and $\deg(ac), \deg(bc), \deg(ad) \leq \deg(D)/2$.

The tabulation of $\ratsfuncq{q}$-isomorphism classes of cubic function fields as performed in \cite{Pieter2,Pieter3,ANTSPie} used the Davenport-Heilbronn bijection between $\ratsfuncq{q}$-isomorphism classes of cubic function fields and a certain collection $\newU$ of $GL_2(\polyringq{q})$-isomorphism classes of binary cubic forms. This set $\newU$ includes all classes of primitive, irreducible binary cubic forms with square-free discriminant, which is all that is required in our context.  The Davenport Heilbronn correspondence simply assigns each form $f(x,y)$ in $\newU$ the irreducible cubic polynomial $f(x,1)$.

In analogy to the number field terminology, a polynomial in \polyringq{q}\ is said to be a \emph{fundamental discriminant} if it is square-free, of degree at least 3, and has leading coefficient 1 or $h$.  In order to compute the $3$-rank of a quadratic function field of square-free discriminant $D$, we need only count the number of  \ratsfuncq{q}-isomorphism classes of cubic function fields of that same discriminant and with at least two infinite places, see Theorem \ref{precisecountoffieldsgiven3rank} below.  
In turn, to list all \ratsfuncq{q}-isomorphism classes of cubic function fields of discriminant $D$, it suffices to enumerate the corresponding unique reduced irreducible binary cubic forms.  To generate all the desired quadratic function fields up to a given discriminant degree bound $B$, we employ the degree bounds stated above (with $\deg(D)$ replaced by $B$) in nested loops over the coefficients of a binary cubic form. For each such form, we check whether it is reduced and has a discriminant $D$ so that $-3D$ is a fundamental imaginary or unusual discriminant of degree at most $B$.  A more precise description of the algorithm is given in Section \ref{sec:alg}.  

Our algorithm for generating quadratic function fields of positive 3-rank relies on key connections between quadratic and cubic fields of the same discriminant.  A modified version of theorem of Hasse \cite{Hasse1} for the function field setting appears in \cite{cuffqi} and gives a precise formula for the number of isomorphism classes of cubic function fields for a fixed square-free discriminant $D$ in terms of the $3$-rank of the quadratic field with discriminant $D$.  Specifically, if $D$ is a square-free polynomial in $\polyringq{q}$ and $K = \ratsfuncqextn{q}{\sqrt{D}}$, then the number of $\ratsfuncq{q}$-isomorphism classes of cubic function fields of discriminant $D$ with at least two infinite places is $(3^{r_3(D)} - 1)/2$, where $r_3(D)$ is the $3$-rank of the ideal class group of the quadratic function field $K$. 




Note that the Hasse count omits the classes of cubic function fields with just one infinite place. To include these classes, we require a more refined count.  Let $n$ be any non-square and suppose that $D$ is unusual. Then the real discriminant $D' = nD$ is said to be the \emph{dual} of $D$, and the real quadratic function field $K' = \ratsfuncqextn{q}{\sqrt{D'}}$ is the \emph{dual} of the unusual quadratic function field $K =  \ratsfuncqextn{q}{\sqrt{D}}$. 
Let $l$ be an odd prime dividing $q+1$. 
If $r$ and $r'$ denote the $l$-rank of $K/\ratsfuncq{q}$ and $K'/\ratsfuncq{q}$, respectively, then $r = r'$ or $r = r' + 1$.  In the latter case, the regulator of $K'/\ratsfuncq{q}$ is a multiple of $l$ (see Lee \cite{LeeY1}).
The cases $r = r' + 1$ and $r = r'$ are referred to as \emph{escalatory} and \emph{non-escalatory}, respectively.

Denote by $(e_1, f_1; \ldots ; e_r, f_r)$ the signature of the place at infinity of \ratsfuncq{q}\ in a finite extension $L$ of $\ratsfuncq{q}$, so $e_i$ is the ramification index and $f_i$ the residue degree of the $i$-th infinite place of $L$ for $1 \leq i \leq r$.  
In the cases of interest, i.e. $-3D$ imaginary or unusual, we have an exact count of the number of $\ratsfuncq{q}$-isomorphism classes of cubic function fields of discriminant $D$ and any given signature.   
The complete statement and proof appear in \cite{cuffqi}.

\begin{theorem} Let $D$ in \polyringq{q}\ be square-free so that $-3D$ is imaginary or unusual. Then the number of $\ratsfuncq{q}$-isomorphism classes of of cubic function fields of discriminant $D$ is \\
$(3^{r_3(-3D)}-1)/2$. Setting $r = r_3(D)$, the possible signatures for these fields and their respective frequencies are as follows:
\begin{itemize}
\item  If $-3D$ is imaginary, then all $(3^r-1)/2$ classes of fields have signature $(1,1;2,1)$.
\item  If $-3D$ is unusual and $q \equiv 1 \pmod{3}$, then all $(3^r-1)/2$ classes of fields have signature $(1,1;1,2)$.
\item If $-3D$ is unusual and $q \equiv -1 \pmod{3}$, then $D$ is the dual discriminant of $-3D$, and there are two possibilities:
\begin{itemize}
\item  In the non-escalatory case, all $(3^r-1)/2$ classes of fields have signature $(1,1;1,1;1,1)$.
\item  In the escalatory case, $(3^r-1)/2$ classes of fields have signature $(1,1;1,1;1,1)$ and the remaining $3^r$ such classes have signature $(1,3)$.
\end{itemize}
\end{itemize}
\label{precisecountoffieldsgiven3rank}
\end{theorem}


\section{The Algorithm and its Complexity}\label{sec:alg}
We now briefly describe our method for tabulating quadratic function fields of imaginary or unusual fundamental discriminant $-3D$ with positive 3-rank up to a given bound $B$ on $\deg(D)$.   The basic algorithm builds on the algorithm for tabulating cubic function fields from \cite{Pieter2,Pieter3,ANTSPie}, except instead of outputting minimal polynomials for all fields, we simply increment a counter for each square-free discriminant found.  The counter and corresponding discriminant values are then output.  
Discriminants and the number of cubic fields with that discriminant are stored in a hash table, and output to a file once the main \emph{for} loops are exited.  Specifically, we loop over each coefficient of a binary cubic form satisfying the bounds given in Section 4.5 of \cite{Pieter2}.  For each binary cubic form $f$ encountered in the loop, we test whether or not $f$ is reduced, is squarefree, $-3D$ is imaginary or unusual, and $\deg(D) \leq B$.  If this is the case, the number of binary cubic forms found with discriminant $D$ is incremented by one in our hash table.   Once the hash table is complete, the value of the counter for each discriminant can then be converted to its proper 3-rank value using Theorem \ref{precisecountoffieldsgiven3rank} depending on which case is appropriate. A modified version of the algorithm, with various improvements including shortening the loop on $d$ as described previously in \cite{Pieter2,Pieter3}, appears in Algorithm \ref{tabimag}.


\begin{algorithm}
\caption{$3$-rank algorithm for computation of quadratic fields where $-3D$ is imaginary (resp. unusual) }
\label{tabimag}
\begin{algorithmic}[1]
\REQUIRE A prime power $q$ not divisible by $2$ or $3$, a primitive root $h$ of $\finfldq{q}$, the set $S=\{1,h,h^2,\ldots h^{(q-3)/2} \}$, and a positive integer $B$.
\ENSURE A table where each entry is a square-free discriminant $D$ and a number of the form \\
$(3^r-1)/2$ with $-3D$ imaginary (resp. unusual), $\sgn(-3D) \in \{1,h\}$ (resp. $\sgn(-3D) = h$), and $\deg(D) = B$.  The positive integer $r$ is the $3$-rank of the quadratic function field \ratsfuncqextn{q}{\sqrt{D}}.
 
\FOR{$\deg(a) \leq B/4$ \textbf{AND} $\sgn(a) \in S$}
\FOR{$\deg(b) \leq B/4$}
\FOR{$\deg(c) \leq B/2 - \max\{\deg(a),\deg(b)$\}}
\STATE $m_1: = 2 (\deg(b) + \deg(c))$;
\STATE $m_2:= \deg(a) + 3 \deg(b)$
\FOR{$i = 0$ to $B/2 - \deg(a)$}
\STATE $m_3: = \deg(a) + \deg(b) + \deg(c) + i$;
\STATE $m_4:= 3 \deg(b) + i$
\STATE $m_5:= 2(\deg(a) + i)$
\STATE $m:= \max \{m_1,m_2,m_3,m_4,m_5\}$
\IF{($m$ is not taken on by a unique term among the $m_i$) \textbf{OR} ($m$ is taken on by a unique term \textbf{AND} $m$ is odd (resp. even) \textbf{AND} $m \leq B$)}
\STATE Compute $P:= b^2 - 3ac$;
\STATE Compute $t_1:= bc$;
\STATE Compute $t_2:= c^2$
\FOR{$\deg(d) = i$}
\STATE Set $f:= (a,b,c,d)$;
\STATE Compute $Q := t_1- 9ad$;
\STATE Compute $R := t_2 - 3bd$;
\STATE Compute $-3D = -3D(f) = Q^2 - 4PR$;
\IF{ $\deg(D) \leq B$ \textbf{AND} 
$-3D$ is imaginary (resp. unusual) 
\textbf{AND}  $-3D$ is fundamental \textbf{AND} $f$ is reduced }
\STATE Increment counter for discriminant $D$ by one in the hash table;
\ENDIF 
\ENDFOR
\ENDIF  
\ENDFOR
\ENDFOR
\ENDFOR
\ENDFOR
\STATE Output hash table;
\end{algorithmic}
\end{algorithm}

The asymptotic complexity of the algorithm for generating fields with positive $3$-rank is the same as for the tabulation algorithm for cubic function fields, namely $O(B^4 q^B)$ field operations \cite{Pieter3}, with the $O$-constant cubic in $q$ when $B$ is odd and quartic in $q$ when $B$ is even. The main difference between the two algorithms is that instead of testing if an equivalence class of binary cubic forms belongs to the Davenport-Heilbronn set, we test if it has square-free discriminant. This requires only one $\gcd$ computation.  The $3$-rank program then stores the discriminant and the $3$-rank data in a hash table, which is output at the end of the algorithm.  

In Belabas \cite{Bela2}, a number of modifications to the basic algorithm for computing the $3$-rank of a quadratic number field were suggested and implemented.  We give a brief summary of these modifications here, and explain why we refrained from making similar changes to our program, but the cost of this is negligible compared to the rest of the algorithm. 

First, we did not use Belabas' ``cluster" approach.  This approach, where one loosens the conditions for a form to be reduced and looks for a large number of forms in a given interval, finally proceeding with class group computations on the clusters found, was not used as we sought to avoid a large number of direct class group computations, except for verification of a small sample of examples.  This does however warrant further investigation.

The other main variants of Belabas' $3$-rank program are dedicated to speeding up the square-free test for integers.  As square-free testing for polynomials is straightforward and efficient, this aspect of Belabas' work was not explored.  

\section{An Improved Algorithm for $-3D$ imaginary}\label{subsec:newalg}

In this section, we will use certain $\ratsfuncq{q}$-automorphisms to speed up our algorithm by a factor of $q$ in most cases.  This requires some additional notation.  Henceforth, let $p$ denote the characteristic of \finfldq{q}.  For any non-constant polynomial $F(t) \in \polyringq{q}$, let $\sgn_2(F)$ denote the coefficient of $t^{\deg(F)-1}$ in $F$ (this is allowed to be zero).   We also require a preliminary lemma, which is easily proved by induction on the degree.

\begin{Lem}
Let $F(t) \in \polyringq{q}$ be a non-zero polynomial. Then $F(t+\beta) = F(t)$ for
all $\beta \in \finfldq{q}$ if and only if $F(t)$ is a polynomial in $t^q-t$.
\label{Translations}
\end{Lem}

\begin{proof}
``Freshmen exponentiation" easily shows that every polynomial in $t^q -  t$ is translation- invariant. The converse certainly holds for constant polynomials.  So let $F(t) \in \polyringq{q}$ be non- constant and assume inductively that the converse statement holds for all polynomials of degree less than $\deg(F)$. The constant coefficient of $F(t)$ is $F (0)$, so obviously $t$ divides $F(t) - F(0)$. Replacing $t$ by $t+\beta$ for any $\beta \in \finfldq{q}$, we see that $t+\beta$ divides $F(t+\beta)-F(0) = F(t)-F(0)$ for all $\beta \in \finfldq{q}$. It follows that

\[t^q - t = \prod_{\beta \in \finfldq{q}} (t - \beta) \]
divides $F(t) - F (0)$. Thus, $F(t) = (t^q - t)G(t) + F (0)$ for some polynomial $G(t) \in \polyringq{q}$ of degree $\deg(F) - q < \deg(F)$. Now again by ``freshmen exponentiation", 
\begin{eqnarray*}
(t^q -t)G(t)+F(0)& = & F(t) =  F(t+\beta) = ((t+\beta)^q -(t+\beta))G(t+\beta)-F(0) \\
& = &(t^q -t)G(t+\beta)-F(0) ,
\end{eqnarray*}
so $G(t + \beta) = G(t)$ for all $\beta \in \finfldq{q}$. By induction hypothesis, $G(t)$ is a polynomial in $t^q - t$, and hence so is $F(t)$.
\end{proof}

The key to our improvements is the following:

\begin{Prop} For every polynomial $F(t) \in \polyringq{q}$ whose degree is coprime to $p$, there exists a unique $\beta \in \finfldq{q}$ such that $\sgn_2(F (t + \beta)) = 0$. 
\label{PolyTransUnique}
\end{Prop}

\begin{proof}
If $d = \deg(F)$, then $\sgn_2(F (t + \beta )) = \sgn_2 (F) + d\beta\,\, \sgn(F)$, which vanishes if and only if $\beta = -\sgn_2 (F)/d \,\sgn(F)$. 
\end{proof}

\begin{Cor}
For every reduced binary cubic form $f=(a,b,c,d)$ over $\polyringq{q}$ with imaginary Hessian and $p \nmid \deg(a)$, there there exists a unique $\beta \in \finfldq{q}$ such that $\sgn_2(a (t + \beta)) = 0$ and the form $f_\beta = (a(t+\beta),b(t+\beta),c(t+\beta),d(t+\beta))$ is reduced. 
\end{Cor}

\begin{proof}
Follows immediately from Proposition \ref{PolyTransUnique} and the fact that translation does not change the conditions on a reduced form as specified in Definition \ref{reduceddefn}.
\end{proof}
Note that if $f$ has discriminant $D(t)$, then $f_\beta$ has discriminant $D(t+\beta)$ which is generally distinct from $D(t)$.  If $f$ is a reduced binary cubic form with unusual Hessian, then translating the coefficients by any $\beta$ in $\finfldq{q}$ produces a partially reduced form, but not necessarily a reduced one.  We were unable to find an efficient way to address this problem.  We will revisit this issue in Section \ref{sec:concl}.

The idea for speeding up the $3$-rank algorithm is to loop only over polynomials $a$ with $\sgn_2(a) = 0$ when $\deg(a)$ is not divisible by $p$. Each such $a$ yields $q$ distinct forms $(a(t + \beta), b(t + \beta), c(t + \beta), d(t + \beta))$ for each non-zero $\beta \in \finfldq{q}$ of respective discriminants $D(t + \beta)$, for which the count is appropriately adjusted afterwards. This can be further improved by using the same idea on the polynomials $d$. That is, the algorithm loops only over pairs $(a, d)$ with $p \nmid \deg(a)$ and $\sgn_2(a) = 0$, or $p \mid \deg(a)$, $p \nmid \deg(d)$ and $\sgn_2(d) \neq 0$, or $p \mid \deg(a)$ and $p \mid \deg(d)$. For the former two, the computational effort decreases by a factor of $q$ if we disregard the computation of all the translates. Only for the last of these three types of pairs is the computational effort the same as in Algorithm \ref{tabimag}. The exact proportion of pairs $(a, d)$ with $p \mid \deg(a)$ and $p \mid \deg(d)$ depends on the residue class of $\lfloor \deg(D)/2  \rfloor \pmod{p}$.  Note that both $a$ and $d$ are non-zero since we only consider irreducible forms.\\
 
For brevity, we introduce the following terminology for cubic forms $f = (a, b, c, d)$. We call the form $f$ a 
\begin{itemize}
\item \emph{type $1$ form} if $p \nmid \deg(a)$ and $\sgn_2 (a) = 0$, or $p \mid \deg(a)$, $p \nmid \deg(d)$ and $\sgn_2 (d) = 0$; 
\item \emph{type 2 form} if $p \mid \deg(a)$ and $p \mid \deg(d)$; 
\item \emph{type 3 form} otherwise, i.e. $p \nmid \deg(a)$ and $\sgn_2(a) \neq 0$, or $p \mid \deg(a)$, $p \nmid \deg(d)$, and $\sgn_2 (d) \neq 0$. \\
\end{itemize}

Note that by Proposition \ref{PolyTransUnique}, if $f = (a, b, c, d)$ is a type 1 form, then for all $\beta \in \finfldq{q}^*$ , $f_\beta = (a(t + \beta), b(t + \beta), c(t + \beta), d(t + \beta))$ is a type 3 form. Moreover, all the $f_\beta$ for $\beta \in \finfldq{q}^*$ are pairwise distinct. Conversely, if $f = (a, b, c, d)$ is a type 3 form, then the forms $f_\beta$ for $\beta \in \finfldq{q}$ are pairwise distinct by Lemma \ref{Translations} (as $a(t)$ and $d(t)$ cannot both be polynomials in $t^q-t$), and by Proposition \ref{PolyTransUnique}, exactly one of the $f_\beta$ is a type 1 form and the others are type 3 forms.  

The revised algorithm, given as Algorithm \ref{tabimagHorner} below, only loops over forms of type 1 and 2 in steps 1--8, and incorporates the discriminant count arising from the type 3  forms via translates in step 11, whereas Algorithm \ref{tabimag} looped over forms of all three types.   Each discriminant $D$ is endowed with two counters. One is the counter keeping track of the number of forms for each discriminant as in step 21 of Algorithm \ref{tabimag}. The other one is a translate counter that counts how often $D$ is encountered as the discriminant of a type 1 form $f = (a, b, c, d)$. If $D(t)$ has translate counter $C_D$ , then each $D(t + \beta)$ with $\beta \in \finfldq{q}^*$ occurs exactly $C_D$ times as the discriminant of the type 3 form $f_\beta = (a(t + \beta), b(t + \beta), c(t + \beta), d(t + \beta))$.
\begin{algorithm}
\caption{Improved $3$-rank algorithm for computation of quadratic fields where $-3D$ is imaginary using Horner's rule for translates}
\label{tabimagHorner}
\begin{algorithmic}[1]
\REQUIRE A prime power $q$ not divisible by $2$ or $3$, a primitive root $h$ of $\finfldq{q}$, the set $S=\{1,h,h^2,\ldots h^{(q-3)/2} \}$, and a positive integer $B$.
\ENSURE A table where each entry is a square-free discriminant $D$ and a number of the form \\
$(3^r-1)/2$ with $-3D$ imaginary, $\sgn(-3D) \in \{1,h\}$, and $\deg(D) = B$.  The positive integer $r$ is the $3$-rank of the quadratic function field \ratsfuncqextn{q}{\sqrt{D}}.
 
\FOR{$\deg(a) \leq B/4$ \textbf{AND} $\sgn(a) \in S$ \textbf{AND} ($\sgn_2(a) = 0$ if $p \nmid \deg(a)$ \textbf{OR} $p \mid \deg(a)$)}
\STATE Execute Steps 2--14 of Algorithm \ref{tabimag}
\FOR{$\deg(d) = i$ \textbf{AND} ($\sgn_2(d) = 0$ if $p \nmid \deg(d)$ \textbf{OR} $p \mid \deg(d)$)}
\STATE Execute Steps 14--19 of Algorithm \ref{tabimag}
\IF{ $\deg(D) \leq B$ \textbf{AND}  
$-3D$ is imaginary (resp. unusual) 
\textbf{AND}  $-3D$ is fundamental \textbf{AND} $f$ is reduced }
\STATE Increment the 3-rank counter for discriminant $D$ by one in the hash table;
\IF{ $p \nmid \deg(a)$ \textbf{OR} $p \nmid \deg(d)$}
\STATE Increment the translate counter of $D(t)$ by 1
\ENDIF 
\ENDIF  
\ENDFOR
\ENDFOR
\FOR{all $D(t)$ in the hash table with translate counter $\geq 1$}
\STATE Compute $D(t+\beta)$ for $\beta \in \finfldq{q}$ using Horner's Rule;
\STATE Increase the counter for $D(t+\beta)$ by the value of the translate counter of $D(t)$;
\ENDFOR
\STATE Output hash table;
\end{algorithmic}
\end{algorithm}
\begin{theorem}
Algorithm \ref{tabimagHorner} is correct.
\end{theorem}

\begin{proof}
We show that Algorithms \ref{tabimag} and \ref{tabimagHorner} have exactly the same output.  Steps 1--8 of Algorithm \ref{tabimagHorner} loop exactly over all the type 1 and 2 forms of degree up to $B$, and no type 3 form. Each type 1 form $f = (a, b, c, d)$ of discriminant $D(t)$ gives rise to $q-1$ type 3 forms $f_\beta = (a(t + \beta), b(t + \beta), c(t + \beta), d(t + \beta))$ of respective discriminants $D(t + \beta)$ for $\beta \in \finfldq{q}^*$ that all generate the same cubic field as $f$. Steps 9-11 generate all these discriminants, and if $D(t)$ has translate counter $C_D$ , i.e.\ is the discriminant of $C_D$ forms $f$ of type 1, then each $D(t + \beta)$ is (in addition to the current value of the counter of $D(t + \beta)$) the discriminant of the $C_D$ forms $f_\beta$ of type 3. So every discriminant output by Algorithm \ref{tabimagHorner} is also output by Algorithm \ref{tabimag}, and with the same $3$-rank count. 

Conversely, let $D$ be a discriminant together with a count that is produced by step 21 of 
Algorithm \ref{tabimag}. Let $D(t)$ be the discriminant of $C_i$ type $i$ forms for $i = 1, 2, 3$. Since steps 1--8 of Algorithm \ref{tabimagHorner} only loop over all type 1 and 2 forms, and no type 3 forms, it produces $D$ with counter $C_1 + C_2$.  Now each of the $C_3$ occurrences of $D(t)$ as the discriminant of a type 3 form $f_i = (a_i , b_i , c_i , d_i )$ corresponds to exactly one occurrence of some translate of $D(t)$ that is the discriminant of a type 1 form as follows. If $p \nmid \deg(a_i)$ and $\sgn_2(a_i) \neq 0$, then there exists a unique $\beta_i \in \finfldq{q}^*$ such that $\sgn_2(a_i (t - \beta_i )) = 0$ by Proposition \ref{PolyTransUnique}. Then $f_i' = (a_i (t - \beta_i ), b_i (t - \beta_i ), c_i (t - \beta_i ), d_i (t - \beta_i ))$ is a type 1 form of discriminant $D(t) = D(t - \beta_i )$. If $p \mid \deg(a_i), p \nmid \deg(d_i )$, and $\sgn_2(d_i) \neq 0$, then there again exists a unique $\beta_i \in \finfldq{q}^*$ such that $\sgn_2 (d_i (t - \beta_i )) = 0$. Then $f_i' = (a_i (t -\beta_i ), b_i (t - \beta_i ), c_i (t- \beta_i ), d_i (t - \beta_i ))$ is again a type 1 form of discriminant $D(t) = D(t -\beta_i )$. Note that these two possibilities are mutually exclusive. For either case, $D(t) = D_i' (t + \beta_i )$ is encountered in steps 9-11. Now each  $D_i'(t)$ occurred exactly $C_{D'}$ times as the discriminant of a type 1 form, and these type 1 forms are exactly the forms $f_i'$. So $D(t) = D_i'(t + \beta_i )$ occurs exactly $C_{D'}$ times as the discriminant of the corresponding type 3 form, and these forms are exactly the forms $f_i$. 
\end{proof}

\section{Complexity of the Improved Algorithm}\label{sec:algcomplexity}

We now analyze the complexity of Algorithm \ref{tabimagHorner}.  As before, let $p$ denote the characteristic of $\finfldq{q}$.  As in \cite{Pieter3}, denote by $\mathcal{F}_s$ the set of binary cubic forms $f = (a, b, c, d)$ over
$\polyringq{q}$ such that $\deg(D(f)) = s$, $\deg(a) \leq s/4$, $\deg(b) \leq s/4$, $\deg(ad) \leq s/2$,
$\deg(bc) \leq s/2$, and $\sgn(a) \in S$. Using the fact that there are $(q-1)/2$ choices for
$\sgn(a)$, $q-1$ choices for $\sgn(d)$, and $q$ choices each for the highest permissable coefficient
of $b$ and $c$, we recall from Lemma 7.1 of \cite{Pieter3} that
\[ \# \mathcal{F}_s = \frac{q^{4 - \delta_s}}{32}
\, s^2 q^s + O(q^s) \]
as $s \rightarrow \infty$, where $\delta_s$ is the parity of $s$. Using this result,
Corollary 7.2 of \cite{Pieter3} established a run time, using the bound $B$ on discriminant degrees, for
the cubic function field tabulation algorithm (Algorithm 2 of \cite{Pieter3}) --- and hence also for Algorithm \ref{tabimag} --- of $O(B^4 q^B)$ operations in $\finfldq{q}$ as $B
\rightarrow \infty$.  Disregarding constants arising from polynomial arithmetic,
the dominant term of the $O$ constant when viewed as a polynomial in $q$ was $q^3/16$ when $B$ is odd and $q^4/32$ when $B$ is even. Our modifications herein improve this constant by a factor of $q$ most of the time; only when $B \equiv 0, 1 \pmod{2p}$ is the speed-up smaller, but still significant.  We begin with an auxiliary lemma.


\begin{Lem} \label{l:aux}
Let $m, n, r \in \mathbb{N}$ with $m \leq n$ and set
\[ N(m,n,r) = \sum_{i \leq m} \, \sum_{i+j \leq n} r^{i+j} \ . \]
Then $N(m,n,r) = \displaystyle \frac{r}{r-1} \, m r^n + O(r^n)$ as $m \rightarrow \infty$.
\end{Lem}
\begin{proof}
\[ N(m,n,r) = \sum_{i=0}^m r^i \, \sum_{j=0}^{n-i} r^j = \sum_{i=0}^m
\frac{r^{n+1} - r^i}{r-1} = \frac{r}{r-1} \, m r^n + O(r^n) \]
as claimed.
\end{proof}

\begin{Cor} \label{c:pairs1}
For $s \geq 0$, let $P(s)$ denote the number of pairs of polynomials $(G, H)$
where $G, H \in \polyringq{q}$, $G$ and $H$ are monic, $\deg(G) \leq s/4$ and $\deg(GH)
\leq s/2$. Then
\[ P(s) = \frac{q}{4(q-1)} \, s q^{\lfloor s/2 \rfloor} + O(q^{s/2}) \]
as $s \rightarrow \infty$.
\end{Cor}
\begin{proof}
The number of monic polynomials in $\polyringq{q}$ of degree at most $i$ is
$q^i$. Hence
\[ P(s) = \sum_{i \leq s/4} \, \sum_{i+j \leq s/2} q^{i+j} = N(s/4, \lfloor s/2
\rfloor, q) \ . \]
The result now follows from Lemma \ref{l:aux}.
\end{proof}

\begin{Cor} \label{c:pairs2}
For $s \geq 0$, let $Q(s)$ denote the number of pairs of polynomials $(G, H)$ where $G, H \in
\polyringq{q}$, $G$ and $H$ are monic, $p \mid \deg(G)$, $p \mid \deg(H)$, $\deg(G) \leq s/4$ and $\deg(GH)
\leq s/2$. Then
\[ Q(s) = \frac{q^p}{4p(q^p-1)} \, s q^{\lfloor s/2 \rfloor - r_s} + O(q^{s/2}) \]
as $s \rightarrow \infty$, where $\lfloor s/2 \rfloor \equiv r_s \pmod{p}$ with $0
\leq r_s \leq p-1$.
\end{Cor}
\begin{proof}
We have
\[ Q(s) = \sum_{pi \leq s/4} \, \sum_{pi+pj \leq s/2} q^{pi+pj}
	= N(s/4p, \lfloor s/2p \rfloor, q^p)
	= \frac{q^p}{4p(q^p-1)} \, s q^{p\lfloor s/2p \rfloor} + O(q^{s/2})\ . \]
It remains to show that $p \lfloor s/2p \rfloor = \lfloor s/2 \rfloor - r_s$ to
deduce the claim from Lemma \ref{l:aux}. Let $\delta_s$ denote the parity of $s$,
so $\lfloor s/2 \rfloor = (s-\delta_s)/2$. Since
\[ \frac{s}{2p} - \frac{s - \delta_s}{2p} = \frac{\delta_s}{2p} \ ,
    \qquad 0 \leq \frac{\delta_s}{2p} < 1 \ , \]
we see that $\lfloor s/2p \rfloor = \lfloor (s-\delta_s)/2p \rfloor$, and hence by definition
of $r_s$,
\[ \left \lfloor \frac{s}{2} \right \rfloor - r_s
	= p \left \lfloor \frac{\lfloor s/2 \rfloor}{p} \right \rfloor
	= p \left \lfloor \frac{(s - \delta_s)/2}{p} \right \rfloor
	= p \left \lfloor \frac{s-\delta_s}{2p} \right \rfloor
	= p \left \lfloor \frac{s}{2p} \right \rfloor \]
as desired.
\end{proof}

\begin{Lem} \label{l:type1}
For any $s$, the number of type 1 forms in $\mathcal{F}_s$ is
\[ N_1(s) = \frac{(q-1)^2 q}{2} \, P(s) \left ( P(s) - Q(s) \right ) \]
where $P(s)$ and $Q(s)$ are defined as in Corollaries \ref{c:pairs1} and
\ref{c:pairs2}, respectively.
\end{Lem}
\begin{proof}
The number of pairs of polynomials $(G, H)$ defined in Corollary \ref{c:pairs1} with $p \nmid
\deg(G)$ or $p \nmid \deg(H)$ is $P(s) - Q(s)$. For those pairs with $p \nmid \deg(G)$, a proportion
of $1/q$ has $\sgn_2(G) = 0$, and for those with $p \mid \deg(G)$ and $p \nmid \deg(H)$, a
proportion of $1/q$ has $\sgn_2(H) = 0$. It follows that the number of $(a, d)$ pairs for which $f =
(a,b,c,d)$ is a type 1 form in $\mathcal{F}_s$ is
\[ \frac{(q-1)^2}{2q} \, \big (P(s) - Q(s) \big ) \ , \]
and the number of $(b,c)$ pairs for such a form is $q^2 P(s)$. Multiplying these two counts yields
the desired result.
\end{proof}

\begin{Lem} \label{l:type2}
For any $s$, the number of type 2 forms in $\mathcal{F}_s$ is
\[ N_2(s) = \frac{(q-1)^2 q^2}{2} \, P(s)Q(s) \]
where $P(s)$ and $Q(s)$ are defined as in Corollaries \ref{c:pairs1} and \ref{c:pairs2}, respectively.
\end{Lem}
\begin{proof}
The number of $(a, d)$ pairs for which $f = (a,b,c,d)$ is a type 2 form in $\mathcal{F}_s$ is \\
$(q-1)^2 Q(s)/2$  and the number of $(b,c)$ pairs for such a form is again $q^2 P(s)$.
\end{proof}

\begin{Cor} \label{c:combined}
The combined number of type 1 and 2 forms in $\mathcal{F}_s$ is
\[ N_1(s) + N_2(s) = \frac{q^{3 - \delta_s}}{32}
    \left ( \frac{(q-1)^2}{pq^{1+r_s}} \frac{q^p}{q^p-1} + 1 \right ) s^2 q^s + O(s q^s) \ . \]
\end{Cor}
\begin{proof}
By Lemmas \ref{l:type1} and \ref{l:type2}, the combined number of type 1 and 2 forms
$\mathcal{F}_s$ is
\[ N_1(s) + N_2(s) = \frac{(q-1)^2 q}{2} \, P(s) \big ( P(s) - Q(s) + qQ(s) \big )
    = \frac{(q-1)^2 q}{2} \, \big ( P(s)^2 + (q-1)P(s)Q(s) \big ) \ . \]
We evaluate each term separately. As before, let $\delta_s$ be the parity of $s$, so that $2 \lfloor
s/2 \rfloor = s - \delta_s$. Then
\[ \frac{(q-1)^2 q}{2} \, P(s)^2 = \frac{(q-1)^2 q}{2} \, \frac{q^2}{16(q-1)^2} \,
        s^2 q^{s - \delta_s} + O(s q^s)
     = \frac{q^{3-\delta_s}}{32} \, s^2 q^{s - \delta_s} + O(s q^s) \ . \]
Setting $(s - \delta_s)/2 \equiv r_s \pmod{p}$ with $0 \leq r_s \leq p-1$ as before, we also obtain
\begin{eqnarray*}
\frac{(q-1)^3 q}{2} \, P(s) Q(s) & = & \frac{(q-1)^3 q}{2} \, \frac{q}{16(q-1)} \,
        \frac{q^p}{4p(q^p-1)} \, s^2 q^{s - \delta_s - r_s} + O(s q^s) \\
    & = & \frac{(q-1)^2 q^{2 - \delta_s - r_s}}{32p} \, \frac{q^p}{q^p - 1} \, s^2 q^s + O(s q^s) \ .
\end{eqnarray*}
Altogether,
\[ N_1(s) + N_2(s) = \frac{q^{3 - \delta_s}}{32}
    \left ( 1 + \frac{(q-1)^2}{pq^{1+r_s}} \frac{q^p}{q^p-1} \right ) s^2 q^s + O(s q^s) \ . \]
\end{proof}

Note that the factor $q^p/(q^p-1)$ is extremely close to 1, even for small values
of $p$ and $q$. For the smallest permissable parameters, when $p = q  = 5$, this
quantity is approximately $1.0003$; for $q = p = 7$, it is roughly $1.000001$.

\begin{theorem} \label{t:runtime}
Assuming standard polynomial arithmetic in $\polyringq{q}$, Algorithm \ref{tabimagHorner}
requires $O(B^4q^B) = O(q^{B+\epsilon})$ operations in $\finfldq{q}$ as $B \rightarrow
\infty$. Here, the $O$-constant is a rational function of $q$ whose dominant term is of
order $q^{4-\delta_B}/p$ if $B \equiv 0, 1 \pmod{2p}$ and of order $q^{3 - \delta_B}$
otherwise, where $\delta_B$ is the parity of $B$.
\end{theorem}
\begin{proof}
The analysis of steps 1--8 of the algorithm proceeds analogously to Corollary 7.3 of \cite{Pieter3}.  These steps run over all the type 1 and type 2 forms in $\mathcal{F}_s$ for $3 \leq s \leq B$. For each such form, the entire collection of polynomial computations requires at most $Ks^2$ field operations for some
constant $K$ that is independent of $B$ and $q$. This holds because all polynomials under
consideration have degree bounded by $s$. 
So the asymptotic run time of steps 1--8 is
\[ T_1(B) = \sum_{s=3}^B \left ( ( N_1(s) + N_2(s)) \cdot K s^2 \right )
	= K \sum_{s=3}^B \left ( C_s s^4 q^s + O(s^3 q^s) \right )
	= K B^4 \sum_{s=3}^B C_s q^s + O(B^3 q^B) \ , \]
where $N_1(s)$ and $N_2(s)$ are given by Lemmas \ref{l:type1} and \ref{l:type2},
respectively, and
\[ C_s = \frac{q^{3 - \delta_s}}{32} \left ( \frac{(q-1)^2}{pq^{1+r_s}} \frac{q^p}{q^p-1} + 1 \right ) \]
by Corollary \ref{c:combined}. Since $\delta_s \geq 0$, $r_s \geq 0$ and $q-1 < q$,
we see that $C_s < C(q)$ where
\[ C(q) = \frac{q^4}{32p} \, \frac{q^p}{q^p-1}  \ . \]
It follows that
\[  \sum_{s=3}^B  C_s q^s < C(q) \sum_{s=0}^{B-2} q^s + C_{B-1} q^{B-1} + C_B q^B
	< \left ( \frac{C(q)}{q(q-1)} + \frac{C_{B-1}}{q} + C_B \right ) q^B \ . \]
The dominant term in $C(q)/q(q-1)$ is $q^2/32p$.

If $B$ is even, then $r_{B-1} \neq 0$, so
\[ \frac{C_{B-1}}{q} + C_B < \frac{q}{32}
	+ \frac{q^3}{32} \left ( \frac{q^{1-r_B}}{p} \frac{q^p}{q^p-1} + 1 \right ) \ . \]
The dominant term in this constant is $q^3(q^{1-r_B}/p + 1)/32$. This is of order $q^4/p$
 if $r_B = 0$, i.e.\ $B \equiv 0 \pmod{2p}$, and $q^3$ otherwise.

If $B$ is odd, then $r_B \neq 0$ and hence $r_{B-1} = r_B - 1$. Hence
\[ \frac{C_{B-1}}{q} + C_B < \frac{q^2}{32} \left ( \frac{q^{2-r_B}}{p} \, \frac{q^p}{q^p-1} + 1 \right )
	+ \frac{q^2}{32} = \frac{q^2}{32} \left ( \frac{q^{2-r_B}}{p} \, \frac{q^p}{q^p-1} + 2 \right ) \ . \]
The dominant term in this constant is $q^2(q^{2-r_B}/p+2)/32$. This is of order $q^3/p$ if $r_B = 1$,
i.e.\ $B \equiv 1 \pmod{2p}$, and $q^2$ otherwise.

Next, we analyze steps 9--11 of Algorithm \ref{tabimagHorner} and show that its run time is negligible compared to that of steps 1--8. For any given degree $s$, the number of discriminants of degree $s$ that steps 9--11 loop over is
certainly bounded above by the number of fundamental imaginary discriminants
of degree $s$, which is $2(q-1)q^{s-1}$, by our remarks at the beginning of Section \ref{ssec:AutoCounts}. This is a very crude upper bound on the number of
$D$ in the table, but it suffices for our purposes.

The run time of step 9--11 is dominated by step 10. Each translate
$D(t+\beta)$ with $\beta \in \mathbb{F}_q^*$ can be computed by applying Horner's Rule as follows: if
$D(t) = a_s t^s + \cdots + a_0$, then initialize $D_0 = a_s$ and compute $D_i = tD_{i-1} + \beta
D_{i-1} + a_{s-i}$ for $1 \leq i \leq s$. Then $D_s = D(t+\beta)$.  For each translate $D(t+\beta)$,
this requires $s$ shifts (whose cost is negligible) and $Ls^2$ operations in $\finfldq{q}$ for some constant
$L$ that is independent of $s$ and $q$. Since there are $q$ choices for $\beta$, the total number of 
field operations required by steps 9--11 is no more than
\[ T_2(B) = \sum_{s=3}^B 2(q-1)q^{s-1} q L s^2 = O(B^2 q^B) \ , \]
which is asymptotically negligible compared $T_1(B)$.
\end{proof}

If $q = p$, then it is most efficient to evaluate $D(t+1)$ from $D(t)$ via
$D_0 = a_s$ and $D_i = tD_{i-1} + D_{i-1} + a_{s-i}$ for $1 \leq i \leq s$.
Applying this technique repeatedly $p-1$ times generates all the polynomials
$D(t+j)$ for $0 \leq j \leq p-1$ using only field additions, no multiplications.

In Tables \ref{t:q!=p} and \ref{t:q=p}, we compare the $O$-constants for the run time of Algorithm \ref{tabimag} as established in Corollary 7.3 of \cite{Pieter3}
with those derived in the proof of Theorem~\ref{t:runtime}. We disregard the constants arising from
polynomial arithmetic (denoted above by $K$ and $L$).

\begin{table}[ht] \caption{Comparison of the O-constants in the run times of Algorithms \ref{tabimag} and \ref{tabimagHorner}, $q \neq p$} \label{t:q!=p}
\begin{center}
\begin{tabular}{|l|c|c|c|c|} \hline
Parity of $B$  & $B \pmod{2p}$                & Algorithm \ref{tabimag}& Algorithm \ref{tabimagHorner} & Speed-up factor \\ \hline\hline
\quad $B$ odd  & $B \not \equiv 1 \pmod{2p}$  & $q^3/16$        & $q^2/16$           & $q$ \\
\quad $B$ odd  & $B \equiv 1 \pmod{2p}$       & $q^3/16$        & $q^3/32p$          & $2p$ \\ \hline
\quad $B$ even & $B \not \equiv 0 \pmod{2p}$  & $q^4/32$        & $q^3/32$           & $q$ \\
\quad $B$ even & $B \equiv 0 \pmod{2p}$       & $q^4/32$        & $q^4/32p$          & $p$ \\ \hline
\end{tabular}
\end{center}
\end{table}

\begin{table}[ht] \caption{Comparison of the O-constants in the run times of Algorithms \ref{tabimag} and \ref{tabimagHorner}, $q = p$} \label{t:q=p}
\begin{center}
\begin{tabular}{|l|c|c|c|c|} \hline
Parity of $B$  & $B \pmod{2p}$                & Algorithm \ref{tabimag}& Algorithm \ref{tabimagHorner}  & Speed-up factor  \\ \hline\hline
\quad $B$ odd  & $B \not \equiv 1 \pmod{2p}$  & $q^3/16$        & $q^2/16$           & $q$ \\
\quad $B$ odd  & $B \equiv 1 \pmod{2p}$       & $q^3/16$        & $3q^2/32$          & $2q/3$ \\ \hline
\quad $B$ even & $B \not \equiv 0 \pmod{2p}$  & $q^4/32$        & $q^3/32$           & $q$ \\
\quad $B$ even & $B \equiv 0 \pmod{2p}$       & $q^4/32$        & $q^3/16$           & $q/2$ \\ \hline
\end{tabular}
\end{center}
\end{table}

\section{Class Group Distribution Results}\label{sec:3rankbg}

We now give some relevant terminology regarding quadratic fields and $3$-ranks, along with some discussion of the Friedman-Washington heuristics and other results on the distribution of class groups of function fields.  We discuss these heuristics, as we will compare the data from our algorithm to these heuristics in order to give some numerical validity to these conjectures in Section \ref{sec:3ranknnumericalresultsimag}.  Before we give an in-depth description of each heuristic/result, we provide a brief overview of the underlying assumptions for each $\ell$-rank result.  Friedman-Washington \cite{FriedWash} assumes $\ell \neq 2$ and $\ell \neq p$, where $p=\charfld(\finfldq{q})$, with imaginary quadratic function fields only.  The result of Ellenberg et al.\ \cite{EllenVenkWest} assumes $q \not\equiv 1 \pmod{\ell}$ and that the extension is imaginary only.  Achter's result \cite{Achter1} assumes $q \equiv 1 \pmod{\ell}$ and that the quadratic function field has only one infinite place (i.e. is either imaginary or unusual).  Garton's result has the same assumptions as that of Ellenberg et al.\ except that $q \equiv 1 \pmod{\ell}$.   

All the aforementioned results apply to Jacobians (i.e. degree zero divisor class groups) of hyperelliptic curves.  It is well-known that the Jacobian and the ideal class group of an imaginary or unusual 
hyperelliptic function field are very closely linked; they are essentially isomorphic
(possibly up to a factor of $\ints/2\ints$), so their respective $\ell$-ranks are equal when $\ell$ is
odd. Therefore, we can use the ideal class group when comparing our data to the heurictics without loss of generality.  We now give an in-depth description of each result.

The Friedman-Washington heuristics \cite{FriedWash} are entirely analogous to the Cohen--Lenstra heuristics.  Loosely speaking, Friedman and Washington predict that given a fixed finite field \finfldq{q} and an abelian $p$-group $H$, where $\ell$ is an odd prime that does not divide $q$, $H$ occurs as the $\ell$-Sylow part of the divisor class group of a quadratic function field over \finfldq{q}\ with frequency inversely proportional to $|\mbox{Aut}(H)|$.  The precise statement is given below.  Denote by $\eta_\infty (\ell)$ the infinite product $\prod_{k \geq 1} (1-\ell^{-k})$.

\begin{conj}[(Friedman and Washington \cite{FriedWash})]
Let $\ell$ be an odd prime that does not divide $q$.
Then a finite abelian group $H$ of $\ell$-power order appears as the $\ell$-Sylow part $Cl_\ell$ of the class group of a quadratic extension $K$ of $\ratsfuncq{q}$ of genus $g_K$ with a frequency inversely proportional to the number of automorphisms of $H$.  That is,
\begin{equation}\lim_{g \rightarrow \infty}\left(\sum_{\substack{K \\ g_K \leq g\\ Cl_\ell \cong H}} 1\left/ \right.\sum_{\substack{K\\g_K \leq g}} 1   \right) = |\mathrm{Aut}(H)|^{-1} \eta_\infty(\ell).\end{equation}
\label{FWconverge}
\end{conj}

A newer, related result, due to Ellenberg, Venkatesh and Westerland \cite{EllenVenk3,EllenVenkWest}, states that the upper and lower densities of imaginary quadratic extensions of $\ratsfuncq{q}$ for which the $\ell$-part of the class group is isomorphic to any given finite abelian $\ell$-group converges to the right-hand side of equation (\ref{FWconverge}), as $q \rightarrow \infty$ with $q \not\equiv 1 \pmod{\ell}$.  

\begin{theorem}[(Ellenberg, Venkatesh and Westerland \cite{EllenVenkWest})]
Let $\ell > 2$ be prime and $A$ a finite abelian $\ell$-group. Write $\delta^+$  (resp.\ $\delta^-$) for the upper density (resp.\ lower density) of imaginary quadratic extensions of \ratsfuncq{q}\ for which the $\ell$-part of the class group is isomorphic to $A$. Then $\delta^+(q)$ and $\delta^-(q)$ converge, as $q \rightarrow \infty$ with $q \not\equiv 0,1 \pmod{\ell}$, to $\eta_\infty (\ell)/|\mathrm{Aut} (A)|$.
\label{Ellenbergetal}
\end{theorem}

The result of Theorem \ref{Ellenbergetal} requires the additional assumption that $q\rightarrow \infty$, but the predicted distribution is what Friedman and Washington assert for fixed $q \not\equiv 1 \pmod{\ell}$.  We note that for fixed 
values of $q$, the Friedman-Washington heuristic is still a conjecture.  Based on Conjecture \ref{FWconverge} and Theorem \ref{Ellenbergetal}, the probability that the $\ell$-rank of an ideal class group of an imaginary quadratic function field is equal to $r$, as given in Cohen and Lenstra \cite{CohLenHeur2} for number fields and in Lee \cite{LeeY2} for function fields, is given by
\begin{equation}
\ell^{-r^2} \eta_\infty (\ell) \prod_{k=1}^{r}(1-\ell^{-k})^{-2}. \label{FWprobfmla}
\end{equation}
For $\ell=3$ and the values of $r=0,1,2,3$ and $4$, we obtain the approximate probabilities $0.56128$, $0.42009$,  $0.019692$, $0.00008739$ and $4.0964 \times 10^{-8}$ respectively.

Achter \cite{Achter1,Achter2} proved a version of the Friedman-Washington heuristic where $q \rightarrow \infty$ inside $\lim_{g \rightarrow \infty}$ with $q \equiv 1 \pmod{\ell}$.   His result is an explicit formula in terms of $p$ for the number 
\begin{equation}\alpha(g, r) = \frac{|\{x \in \mathrm{Sp}_{2g}(\mathbb{F}): \, \ker(x-id) \simeq \mathbb{F}^r\}|}{|\mathrm{Sp}_{2g}(\mathbb{F})|},\label{AchterFmla}\end{equation}
where $\mathrm{Sp}_{2g}(\mathbb{F})$ denotes the group of $2g$ by $2g$ symplectic matrices over a field $\mathbb{F}$ or order $\ell$.  The function $\alpha(g,r)$ expresses similar probabilities as in the original Friedman-Washington heuristic in the case $q \equiv 1 \pmod{\ell}$.
Achter's result \cite{Achter1} proves a weaker version of Friedman and Washington's original conjecture since it requires $q \rightarrow \infty$ in addition to $g\rightarrow \infty$.

Empirical evidence (see Malle \cite{Malle,Malle2}) strongly suggests that the presence of $\ell$-th roots of unity in the base field skews the distribution of $\ell$-rank values of quadratic number fields.    In order to account for this discrepancy for number fields, Malle \cite{Malle2} proposed alternative conjectural formulas to cover the case where primitive $\ell$-th roots of unity lie in the base field, provided primitive $\ell^2$-th roots of unity that are themselves not also $\ell$-th roots of unity do not lie in the base field.  

A different probability distribution, one consistent with Malle's conjectures \cite{Malle,Malle2} and Achter's results \cite{Achter1,Achter2}, is obtained for $q \equiv 1 \pmod{\ell}$ due to the presence of $\ell$-th roots of unity.  
The analogous weak Friedman-Washington result in this case is due to Garton \cite{Garton,Garton2}.  Garton's result, again with $q \rightarrow \infty$ as well as $g \rightarrow \infty$, predicts that the upper and lower densities of imaginary quadratic extensions of $\ratsfuncq{q}$ for which the $\ell$-part of the class group is isomorphic to any given finite abelian $\ell$-group converges to Malle's formula for number fields, as $q \rightarrow \infty$ with $q \equiv 1 \pmod{\ell}$.  This is summarized in the following theorem.

\begin{theorem}[(Garton \cite{Garton,Garton2})]
Let $\ell > 2$ be prime with $\ell \nmid q$ and $A$ a finite abelian $\ell$-group. Write $\delta^+$  (resp.\ $\delta^-$) for the upper density (resp.\ lower density) of imaginary quadratic extensions of \ratsfuncq{q}\ for which the $\ell$-part of the class group is isomorphic to $A$. Then $\delta^+(q)$ and $\delta^-(q)$ converge, as $q \rightarrow \infty$ with $q \equiv 1 \pmod{\ell}$, to
 
 \[c \cdot \left( \prod_{j \geq1}^{r} \frac{\ell^{j}-1}{\ell^r}\cdot \frac{1}{|H||\mathrm{Aut} (A)|}\right), 
\mbox{ where } c = \left(\prod_{j \geq1} (\ell^{-j}+1)\right)^{-1}.\]
\label{GartonDist}
\end{theorem}

From the distribution given in Theorem \ref{GartonDist}, the probability that a quadratic function field has $\ell$-rank equal to $r$ in the case $q \equiv 1 \pmod{\ell}$ is given by
\begin{equation}
\ell^{-(r^2+r)/2} \frac{\eta_\infty (\ell)}{\eta_\infty(\ell^2)}\prod_{k=1}^{r}(1-\ell^{-k})^{-1}. \label{Malleprobfmla}
\end{equation}
For $\ell=3$ and the values of $r=0,1,2,3$ and $4$, we obtain the approximate probabilities $0.64032
$, $0.31950$,  $0.03994$, $1.5361 \times 10^{-3}$ and $1.9201 \times 10^{-5}$ respectively.

We note that Achter's function $\alpha$ converges to a formula of Malle for number fields as $g\rightarrow \infty$, giving additional evidence for a stronger version of Garton's result, where the dependence on $q \rightarrow \infty$ is removed.
This result is Proposition 3.1 of Malle \cite{Malle2}.    In Section \ref{sec:3ranknnumericalresultsimag}, 
we will use the numerical data we generated to see how well the data matches the the Friedman-Washington/Ellenberg-Venkatesh-Westerland result for $q \equiv -1 \pmod{\ell}$ and the Garton distribution result for $q \equiv 1 \pmod{\ell}$ when $q$ is fixed.

\section{Numerical Results}\label{sec:3ranknnumericalresultsimag}

Tables \ref{3rankimagtable} and  \ref{3rankunustable} present the results of our computations for
the $3$-rank counts of quadratic function fields for $q=5,7,11,13$ using Algorithm \ref{tabimagHorner} for imaginary discriminants and Algorithm \ref{tabimag} for unusual discriminants $-3D$. 
We implemented our
counting algorithm using the C++ programming language coupled
with the number theory library NTL \cite{Shoup1}.  The lists of quadratic fields and their (positive) $3$-rank values were computed on a 32 core 2.0 GHz Intel Xeon X7550 machine running Unix with 256 GB of RAM. 
Each table entry consists of the base field size $q$, the degree bound on the discriminant $D$ and the corresponding genus $g$, the total number of square-free discriminants of that degree, the $3$-rank of $\ratsfuncqextn{q}{\sqrt{D}}$, the total number of fields with given $\deg(D)$ value found with that $3$-rank, and the total elapsed time to find all quadratic function fields with the given degree and $3$-rank at least $1$.  

\begin{center}
\begin{longtable}{|c|c|c|c|c|c|} 
\caption{$3$-ranks of quadratic function fields over $\finfldq{q}$ with $-3D$ imaginary}\label{3rankimagtable} \\

\hline \multicolumn{1}{|c|}{\,\,$q$\,\,} & \multicolumn{1}{|c|}{$\deg(D)$, $g$} & \multicolumn{1}{|c|}{\# of $D$} & \multicolumn{1}{c|}{$3$-rank} & \multicolumn{1}{c|}{Total} & 
\multicolumn{1}{c|}{Total elapsed time} \\ \hline \endfirsthead

\multicolumn{6}{c}%
{{\bfseries \tablename\ \thetable{} -- continued from previous page}} \\
\hline \multicolumn{1}{|c|}{$q$} & \multicolumn{1}{|c|}{$\deg(D)$, $g$} & 
\multicolumn{1}{c|}{\# of $D$} &
 \multicolumn{1}{c|}{$3$-rank} & \multicolumn{1}{c|}{Total} & 
\multicolumn{1}{c|}{Total elapsed time} \\ \hline
\endhead

\hline \multicolumn{6}{|c|}{{Continued on next page}} \\ \hline
\endfoot

\hline 
\endlastfoot
\,\,$5$\,\, & $3$, $1$& $200$ & $1$ & $80$ & 0.00 seconds \\ \cline{2-6}
&$5$, $2$& $5000$ &$1$ & $1600$ & 0.18 seconds  \\ \cline{4-5}
&                &  & $2$ & $10$      &                       \\ \cline{2-6}
&$7$, $3$&125000 &$1$ &  $46840$   & 5.45 seconds \\ \cline{4-5}
&                & & $2$ & $1180$      &             \\ \cline{2-6}
&$9$, $4$& 3125000& $1$ &  $1297120$  &  3 minutes, 22 sec\\ \cline{4-5}
&                 & & $2$ & $51300$       &   \\ \cline{4-5}
&                 & & $3$ & $40$              & \\  \cline{2-6}
&$11$, $5$& 78125000& $1$ & $31730080$  & 2 hours,\\ \cline{4-5}
&                   & & $2$ & $1167200$   &  15 min, 26.21 sec \\ \cline{4-5}
&                   & & $3$ & $1880$          &  \\ \cline{2-6}
&$13$, $6$&1953125000 & $1$ & 806759000& 3 days, 23 hours,\\ \cline{4-5}
&                   & & $2$ & 33601470  &  40 min \\ \cline{4-5}
&                   & & $3$ &  88680         & 
   \\ \hline
\,\,$7$\,\,  & $3$, $1$& 588& $1$ & $196$  & 0.04 seconds\\ \cline{4-5}
&                  & & $2$ & $14$ & \\ \cline{2-6}
&$5$, $2$& 28812 &$1$ & $8400$ & 2.48 seconds  \\ \cline{4-5}
&                  & & $2$ & $588$     &                      \\\cline{2-6}
&$7$, $3$& 1411788 & $1$ &  $432348$ & 2 minutes, 20 sec   \\  \cline{4-5}
&                 & & $2$ & $42924$      &  \\ \cline{4-5}
&                 & & $3$ & $840$           & \\ \cline{2-6}
&$9$, $4$& 69177612 &$1$ &  $21996478$  &  2 hours,\\ \cline{4-5}
&                  & & $2$ & $2675302$   &  44 min, 24 sec  \\ \cline{4-5}
&                  & & $3$ & $90874$        & \\ \cline{4-5}
&                  & & $4$ &  $588$      & \\ \cline{2-6}
&$11$, $5$& 3389702988  &$1$ & 1072738086    & 8 days, 14 hours,\\ \cline{4-5}
&                  & & $2$ & 126751170  &   48 min, 9 sec \\ \cline{4-5}
&                  & & $3$ &       4078662  & \\ \cline{4-5}
&                  & & $4$ &      27174       & \\ \hline \hline
\,\,$11$\,\,&$3$, $1$& 2420& $1$ & $1100$  & 0.51 seconds\\ \cline{2-6}
&$5$, $2$& 292820&  $1$ & $110000$  &  1 minutes, 12 sec    \\   \cline{4-5}                                                                                                           
&                 & & $2$ & $2970$     &                    \\ \cline{2-6}
&$7$, $3$& 35431220 &$1$ &  14186480 & 2 hours,   \\  \cline{4-5}
&                & & $2$ &   506220    &  40 min, 4 sec 
\\ \cline{4-5}
&                & & $3$ &  660   & \\ \cline{2-6}
&$9$, $4$& 4287177620 &$1$ &   1796730320 &   24 days, 23 hours
 \\ \cline{4-5}
&                  & & $2$ &   81402640 &    \\ \cline{4-5}
&                  & & $3$ &    288200     & \\ \hline \hline
\,\,$13$\,\,& $3$, $1$& 4056& $1$ & $1352$  & 1.35 seconds\\ \cline{4-5}
&                & & $2$ & $130$   & \\ \cline{2-6}
&$5$, $2$& 685464&$1$ & $209352$  &  4 minutes, 25 sec \\ \cline{4-5}
&                &  & $2$ & $20046$      &  \\ \cline{4-5}
&                 & & $3$ & $312$       &        \\ \cline{2-6}
&$7$, $3$ & 115843416& $1$ &  36281076  & 14 hours  \\  \cline{4-5}
&                & & $2$ &  4009330       &     38 min, 34 sec
\\ \cline{4-5}
&                & & $3$ & 108108 & \\  \cline{4-5}
&                & & $4$ & 494   & \\ \cline{2-6}
\end{longtable}
\end{center}

\begin{center}
\begin{longtable}{|c|c|c|c|c|c|} 
\caption{$3$-ranks of quadratic function fields over $\finfldq{q}$ with $-3D$ unusual}\label{3rankunustable} \\

\hline \multicolumn{1}{|c|}{\,\,$q$\,\,} & \multicolumn{1}{|c|}{$\deg(D)$, $g$} & \multicolumn{1}{|c|}{\# of $D$} & \multicolumn{1}{c|}{$3$-rank} & \multicolumn{1}{c|}{Total} & 
\multicolumn{1}{c|}{Total elapsed time} \\ \hline \endfirsthead

\multicolumn{6}{c}%
{{\bfseries \tablename\ \thetable{} -- continued from previous page}} \\
\hline \multicolumn{1}{|c|}{$q$} & \multicolumn{1}{|c|}{$\deg(D)$, $g$} & 
\multicolumn{1}{c|}{\# of $D$} &
 \multicolumn{1}{c|}{$3$-rank} & \multicolumn{1}{c|}{Total} & 
\multicolumn{1}{c|}{Total elapsed time} \\ \hline
\endhead

\hline \multicolumn{6}{|c|}{{Continued on next page}} \\ \hline
\endfoot

\hline 
\endlastfoot
\,\,$5$\,\,& $4$, $1$& $500$ & $1$ & $200$ & 0.54 seconds \\ \cline{2-6}
&$6$, $2$& 12500 &$1$ & $4780$  & 12.17 seconds \\ \cline{4-5}  
&                &  & $2$ & $100$   &                      \\   \cline{2-6}
&$8$, $3$& 312500 & $1$ &  $115460$  &  7 minutes, 48.69 sec\\  \cline{4-5}
&                & & $2$ & $2205$      &          \\ \cline{2-6}
&$10$, $4$& 7812500 & 1 &  $3240340$  &  3 hours,  \\ \cline{4-5}
&                  &  & 2 & $128160$  & 9 min, 18.06 sec   \\ \cline{4-5}
&                  &  & 3 & $100$ & \\ \hline \hline
\,\,$7$\,\,  & $4$, $1$& 2058 & $1$ & $672$ & 7.08 seconds \\ \cline{4-5}
&                  & & $2$ & $42$ & \\ \cline{2-6}
&$6$, $2$& 100842 &$1$ & $30989$  & 5 minutes, 42.99 sec   \\ \cline{4-5}
&                 & & $2$ & $3115$     &    \\  \cline{4-5}
&                 & & $3$ &  $63$  & \\ \cline{2-6}
&$8$, $3$& 4941258 &$1$ & $1510026$  & 6 hours, \\ \cline{4-5}
&                & &  $2$ &  $142632$     &  36 min, 20.94 sec  \\ \cline{4-5}
&                & &   $3$  & $2310$        &  \\ \hline \hline
\end{longtable}
\end{center}

Timings for Algorithm \ref{tabimagHorner} are compared to those of Algorithm \ref{tabimag} in Table \ref{ModvsBasicimag} for $-3D$ imaginary with $q=5,7,11,13$.   As seen in the third column of these tables, the modified algorithm is a significant improvement over the basic algorithm.  These improvements get better as $q$ increases.  This is expected as the the complexity analysis for Algorithm \ref{tabimagHorner} predicts an improvement in timings by roughly a factor of $q$, with an improvement of $2q/3$ for degree 11 fields over $\finfldq{5}$.  The actual speed up is below the factor predicted by Table \ref{t:q=p}, but this is likely due to the fact that our degree bounds are too small for the asymptotics to take effect, so the error terms have a significant effect on the run times. \\

\begin{center}
\begin{longtable}{|c|c|c|c|c|} 
\caption{Timings of Algorithm \ref{tabimag} versus Algorithm \ref{tabimagHorner}}\label{ModvsBasicimag} \\

\hline \multicolumn{1}{|c|}{$q$} & \multicolumn{1}{|c|}{Deg.\ bd.} & \multicolumn{1}{|c|}{Alg.\ \ref{tabimag} times} & \multicolumn{1}{c|}{Alg.\ \ref{tabimagHorner} times} & 
\multicolumn{1}{c|}{Alg.\ \ref{tabimag}/Alg.\ \ref{tabimagHorner}} \\ \hline \endfirsthead

\multicolumn{5}{c}%
{{\bfseries \tablename\ \thetable{} -- continued from previous page}} \\
\hline \multicolumn{1}{|c|}{$q$} & \multicolumn{1}{|c|}{Degree bd.} & \multicolumn{1}{|c|}{Algorithm \ref{tabimag} times} & \multicolumn{1}{c|}{Algorithm \ref{tabimagHorner} times} & 
\multicolumn{1}{c|}{Alg.\ \ref{tabimag}/Alg.\ \ref{tabimagHorner}} \\ \hline
\endhead

\hline \multicolumn{5}{|c|}{{Continued on next page}} \\ \hline
\endfoot

\hline 
\endlastfoot
$5$ & $3$  & 0.01 seconds & 0.00 seconds & -- \\ \cline{2-5}
&$5$  & 0.59 seconds & 0.18 seconds & 3.28 \\ \cline{2-5}
&$7$  & 17.33 seconds & 5.45 seconds & 2.99 \\ \cline{2-5}
&$9$  &  11 minutes, 59 sec & 3 minutes, 22 sec & 3.57\\ \cline{2-5}
&$11$  & 5 hours,   & 2 hours, & 2.25\\ 
 &           & 5 min, 19 sec   & 15 min, 26 sec & \\  \cline{2-5}

 &$13$  & ---   & About 4 days & ---\\ \hline

 $7$ & $3$  & 0.09 seconds  & 0.04 seconds & 2.25 \\ \cline{2-5}
& $5$  & 10.75 seconds & 2.48 seconds & 4.33 \\\cline{2-5}
& $7$  & 10 minutes, 2 sec & 2 minutes, 20 sec & 4.30 \\ \cline{2-5}  
& $9$  & 13 hours,  & 2 hours,& 5.03 \\  
&          &   47 min, 6 sec   &  44 min, 24 sec  &   \\ \cline{2-5}
&$11$  & ---   & 8 days, 15 hours & ---\\ \hline

$11$ & $3$  & 1.34 seconds  & 0.51 seconds & 2.63 \\ \cline{2-5}
& $5$  & 6 minutes, 43 sec & 1 minute, 12 sec & 5.58 \\\cline{2-5}
& $7$  & 17 hours, & 2 hours,  & 6.39\\
&          &   2 min, 52 sec    & 40 min, 4 sec  &   \\  \cline{2-5}
&$9$  & ---   & 24 days, 23 hours & ---\\ \hline
$13$ & $3$  & 4.26 seconds  & 1.35 seconds & 3.16 \\ \cline{2-5}
& $5$  & 25 minutes, 44 sec & 4 minutes, 25 sec & 5.81 \\\cline{2-5}
& $7$  & 3 days, 14 hours, & 14 hours, & 5.90 \\ 
&         &     24 min, 22 sec      &    38 min, 34 sec    &    \\ \cline{2-5}
\end{longtable}
\end{center}

We were able to produce examples of escalatory and non-escalatory cases over \finfldq{5}; this was ascertained by computing the class groups of these particular quadratic function fields and their corresponding dual discriminants.  For example, the real quadratic function field of discriminant $D = t^{10} + 2t^9 + t^8 + 4t^7 + 2t^6 + 3t^5 + 3t^4 + 4t^3 + 3t^2 + t$ has $3$-rank $2$, but its unusual dual field of discriminant  $-3D$ has $3$-rank $3$.  

For each of the finite fields specified above, $h=2$ was chosen as a primitive root of $\finfldq{q}$ with the exception of $\finfldq{7}$, where $h=3$ was chosen.  This completely determines the set $S$ specified in Section \ref{sec:prel}.  These sets were $\{1,2\}$, $\{1,2,3\}$, $\{1,2,4,5,8\}$ and $\{1,2,3,4,6,8\}$ for $\finfldq{5}$, $\finfldq{7}$, $\finfldq{11}$ and $\finfldq{13}$ respectively.

We note that our algorithm in both the imaginary and the unusual case is particularly successful at finding quadratic function fields of high $3$-rank and small genus because our method is exhaustive.  This means that any examples of a given $3$-rank value found by our algorithm are minimal, in the sense that any quadratic field with the same $3$-rank must have genus at least as large as the examples found by our algorithm.  For example, our algorithm beats the Diaz~y~Diaz method, (\cite{Diaz1}; adapted to quadratic function fields in \cite{BauJac}), in the sense that fields of $3$-rank equal to $3$ were found for genus $4$ fields over $\finfldq{5}$ using our algorithm.  The minimal genus yielded by the Diaz~y~Diaz method in \cite{BauJac} for quadratic function fields over \finfldq{5}\ with $3$-rank equal to $3$ is $g=5$.  For $\finfldq{11}$ and $\finfldq{13}$ and $3$-rank values $3$ and $4$, we also found examples of smaller genus than those given in \cite{BauJac}.   The minimal genus values were the same for $\finfldq{7}$ and $3$-rank $4$ for both methods.   

An explicit comparison of our method to Diaz~y~Diaz's algorithm is given in Tables \ref{DiazComp} and \ref{DiazCompUnus}.  The third column denotes the minimal genus found with Diaz~y~Diaz' method yielding the corresponding 3-rank specified in column 2.  The fourth column denotes the minimal genus found with our method for the same $3$-rank.  The fifth column gives an example discriminant of minimal degree with given $3$-rank found by our method.  Table \ref{DiazCompUnus} treats the case where $-3D$ is unusual.
These tables indicate a genuine improvement over previous methods in this regard with respect to finding minimal genera with high $3$-rank.  

\begin{table}
\caption{Minimal genera:  our method versus Diaz y Diaz's method ($-3D$ imaginary)} 
\centering
\begin{tabular}{|c|c|c|c|r|}\hline
\multicolumn{1}{|c|}{$q$} & \multicolumn{1}{c|}{$3$-rank} & \multicolumn{1}{c|}{min $g$ (DyD)} & \multicolumn{1}{c|}{min $g$ (RJS)} & \multicolumn{1}{c|}{Example $D(t)$}\\ \hline

$5$& $3$  &  $5$ & $4$    &   $t^9 + 2t^6 + 2t^3 + 3$  \\ \hline
$7$ & $3$ &  $4$ & $3$  &  $3t^7 + 4t^6 + t^2 + 5t + 1$  \\ \cline{2-5}
        & $4$ & $4$ &  $4$ & $3t^9 + 5t^8 + 2t^7 +  $ \\
        &         &          &         & $2t^4 + 3t^3 + 5t^2 + t + 1$ \\ \hline
$11$ & $3$ &  $4$  & $3$  &  $2t^7 + 3t^6 + t^5 + 8t^4 +$ \\ 
           &      &         &          & $3t^3 + t^2 + 3t + 8$ \\ \hline
$13$ & $3$ &  $4$ &  $2$  &  $t^5 + 10t^3 + 8t^2 + 1$ \\ \cline{2-5}
          &   $4$   &  $4$  & $3$  &  $t^7 + 10t^6 + 12t^5 + 5t^4 +$ \\
  & & & & $3t^3 + 11t^2 + 4t + 10$  \\ \hline
\end{tabular}
\vspace{3mm}
\label{DiazComp}
\end{table}

\begin{table}
\caption{Minimal genera:  our method versus Diaz y Diaz's method ($-3D$ unusual)} 
\centering
\begin{tabular}{|c|c|c|c|r|}\hline
\multicolumn{1}{|c|}{$q$} & \multicolumn{1}{c|}{$3$-rank} & \multicolumn{1}{c|}{min $g$ (DyD)} & \multicolumn{1}{c|}{min $g$ (RJS)} & \multicolumn{1}{c|}{Example $D(t)$}\\ \hline

$5$& $3$ & $5$ & $4$  &   $2t^{10} + 4t^9 + 2t^8 + 3t^7 + 4t^6 +$    \\ 
            &         &     &                        &   $t^5 + t^4 + 3t^3 + t^2 + 2t$   \\ \hline
$7$ & $3$ &  $4$ & $3$   &    $3t^8 + t^6 + 5t^2 + 4t + 6$ \\ \hline
\end{tabular}
\vspace{3mm}
\label{DiazCompUnus}
\end{table}

\section{Observations on our Numerical Data}\label{ssec:AutoCounts}
\subsection{Divisibility Properties of the Field Counts}

We note that various values in Table \ref{3rankimagtable} have some interesting divisibility properties.  For instance, all entries in column 5 are divisible by $q$, entries in column 5 for odd degree discriminants are divisible by 2, and some entries in column 5 are divisible by $q-1$.  These properties can be explained via automorphisms in the field $\ratsfuncq{q}$. 
We now briefly discuss the effect of $\ratsfuncq{q}$-automorphisms on the number of discriminants in the various columns of Table \ref{3rankimagtable}, thereby solving an open question from \cite{Pieter2,Pieter3}.  
We note that the counts in column 3 are always of the form $k(q - 1)q^{d-1}$ where $d$ is the discriminant degree and $k = 2$ for $d$ odd and $k = 1$ for $d$ even.  
A standard inclusion-exclusion argument proves that the number of monic square-free polynomials of degree $d \in \mathbb{N}$ with coefficients in $\finfldq{q}$ is $(q - 1)q^{d-1}$; see Exercise 2 of Section 4.6.2, p.\ 456, and its solution on p. 679, of \cite{Knuth2}.  
The value of $k$ is also easily explained: for odd degree, there are two choices for $\sgn(D)$, namely 1 and $h$, whereas for $d$ even, there is one, namely $h$.  Note also that every quadratic function field over $\finfldq{q}$ is the function field of a curve of the form $y^2 = D(t)$ where $D(t)$ is a fundamental discriminant.   Our first result explains why the our 3-rank totals in column 5 of Table \ref{3rankimagtable} are divisible by $q$.

\begin{Prop}
Assume that $d \geq 3$ and $q$ does not divide $d$. Then for any $3$-rank $r$, the number of
imaginary/unusual/real fundamental discriminants of degree $d$ that define a quadratic function field of $3$-rank $r$ is a multiple of $q$.
\label{divisbyq}
\end{Prop}

\begin{proof}
Let $D(t) \in \polyringq{q}$ be a fundamental discriminant of degree $d$. Then the polynomials $D(t+\beta)$ with $\beta \in \finfldq{q}$ are all fundamental discriminants of the same type (imaginary/unusual/real). Since $q \nmid d$, $D(t)$ cannot be a polynomial in $t^q - t$. By Lemma \ref{Translations}, the polynomials $D(t+\beta)$ with $\beta \in \finfldq{q}$ are pairwise distinct. Furthermore, the $q$ curves $y^2 = D(t+\beta)$ are all isomorphic and hence define the same quadratic function field.
\end{proof}

It remains to explain why the entries of column 5 for which the discriminant degree is a multiple of $q$ are also divisible by $q$. Note that the proof of Proposition \ref{divisbyq} remains valid for fundamental discriminants $D$ such that $q \mid \deg(D)$ but $D$ is not a polynomial in $t^q - t$. For the remaining discriminants, the divisibility by $q$ is simply an artifact of the small parameter sizes and no longer holds for larger values of $q$ and $\deg(D)$. For example, the quadratic function field $\ratsfuncqextn{q}{\sqrt{t^q - t}}$ has 3-rank 3 for $q = 11$ and 3-rank 6 for $q = 13$.   Explicit formulas for the class numbers of $\ratsfuncqextn{q}{\sqrt{t^q-t+j}}$ with $q$ an odd prime can be found in \cite{Duursma1}.

Next, we explain why each entry in column 5 of Table \ref{3rankimagtable} is even.

\begin{Prop}
Assume that $q$ is odd and let $d \geq 3$ be odd. Then for any $3$-rank $r$, the
number of imaginary fundamental discriminants of degree $d$ that define a function field of $3$-rank $r$
is even.
\label{divisby2}
\end{Prop}

\begin{proof}
Suppose $D(t)$ is monic. Then the curves $y^2 = D(t)$ and $(h^{(d-1)/2}y)^2 = D(ht)$ are isomorphic
over $\polyringq{q}$ and hence have the same quadratic function field. The result now follows since the distinct polynomials $D(t)$ and $h^{1-d}D(ht)$ are imaginary fundamental discriminants with respective
leading coefficients 1 and $h$.
\end{proof}

We remark that many, but not all, of the entries in column 5 of Table \ref{3rankimagtable} are divisible by $q - 1$. This is likely due to the fact that most of the time, for the corresponding discriminants $D(t)$, the
polynomials $D(at)$ with $a \in \finfldq{q}^*$ are pairwise distinct.
For example, we have the following result:

\begin{Prop}
Assume that $q$ is an odd prime and $D(t) \in \polyringq{q}$ a fundamental discriminant with $D(t)$ monic if $\deg(D)$ is odd. Suppose that for at least one $j \in \{0, 1, \ldots d - 1\}$ with $d - j$ coprime to $q - 1$, the coefficient of $t^j$ in $D(t)$ is non-zero. Then the $q - 1$ hyperelliptic curves $y^2 = h^{-2 \left\lfloor id/2 \right\rfloor}D(h^it)$ for $0 \leq i \leq q-2$ are pairwise distinct and isomorphic, and hence have the same quadratic function field.
\label{divisbyqminus1}
\end{Prop}

\begin{proof}
For brevity, set $F_i(t) = h^{-2\left \lfloor id/2  \right\rfloor}D(h^it)$. If $d$ is even, then all the $F_i(t)$ have the same leading coefficient, namely $\sgn(D)$. So they are either all real or all imaginary fundamental discriminants.  If $d$ is odd, then $F_i(t)$ is monic if $i$ is even and has leading coefficient $h$ if $i$ is odd; in either case, all $F_i(t)$ are imaginary fundamental discriminants. Hence, the hyperelliptic curves 
$y^2 = F_i(t)$ are well-defined and isomorphic over $\finfldq{q}$, and hence they all have the same function field.

It remains to show that the polynomials $F_i(t)$ are pairwise distinct. If $a_j$ is the coefficient of $t^j$ in $D(t)$, then the coefficient of $t^j$ in $F_i(t)$ is $a_jh^{ij-2\left\lfloor id/2 \right\rfloor}$. So assume that $a_j \neq 0$ and $\gcd(q-1, d-j) = 1$, and suppose that $F_i(t) = F_k(t)$ with $0 \leq i,k \leq q - 2$. Then
\begin{equation}
h^{ij-2\left\lfloor id/2 \right\rfloor} = h^{kj-2\left\lfloor kd/2 \right\rfloor}.
\label{floororder}
\end{equation}
If $d$ is even, then (\ref{floororder}) reduces to $h^{(j-d)(i-k)} = 1$, so $(j - d)(i - k) \equiv 0 \pmod{q-1}$. Since $j - d$ is coprime to $q - 1$, this forces $i = k$. 

If $d$ is odd, then $i$ and $k$ must be of the same parity, as otherwise one of $F_i(t)$ would be monic and the other would have sign equal to $h$. The exact same argument as for $d$ even again forces $i = k$.
\end{proof}

Most discriminants satisfy this condition of Proposition \ref{divisbyqminus1}. If those discriminants that do not satisfy this condition correspond to function fields of 3-rank $0$ over $\finfldq{q}$, then the number of fields with positive 3-rank is divisible by $q - 1$ for a given $D$.  Note, however, that this is not always the case; counterexamples include $q = 5, \deg(D) \in \{5, 8\}$ and $q = 7, \deg(D) \in \{3, 6\}$. 

\subsection{Comparison of the 3-rank Data to Heuristics}

\begin{figure}
\caption{Comparison of actual 3-rank data with predicted value (\ref{Ellenbergetal}) over $\finfldq{5}$}
\centering
\includegraphics[width=350pt, height=233pt]{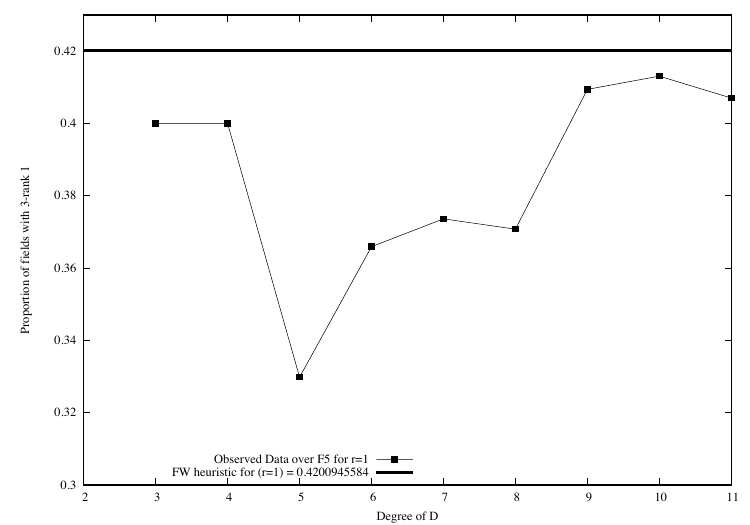}
\label{F5FriedWashComp}
\end{figure}

Since our algorithm is exhaustive for imaginary or unusual discriminants $-3D$ of fixed degree, we compared our data to the distribution suggested by the Friedman-Washington heuristic (Conjecture \ref{FWconverge}) and the Ellenberg-Venkatesh-Westerland result (Theorem \ref{Ellenbergetal}), in an effort to provide numerical evidence in support of the results' validity in the case of fixed $q$.  
A comparison of our $3$-rank data for imaginary and unusual discriminants for the field $\finfldq{5}$ to the expected value predicted by the aforementioned results is given in Figure~\ref{F5FriedWashComp}.   We compared our data to the value given by (\ref{FWprobfmla}) in the imaginary and unusual cases for $p=3$.
The solid line denotes the value of (\ref{FWprobfmla}) for $r=1$ and the dotted line denotes the proportion of $3$-rank one values found up to a given bound found by our algorithm.
As seen in Figure~\ref{F5FriedWashComp}, our data mostly agrees with the predicted asymptotic probability (\ref{FWprobfmla}).  The data for $\finfldq{11}$ (figure omitted)
does not agree as closely, but this is likely because computations were not carried out far enough to obtain a sufficient sample size.   

\begin{figure}
\caption{Comparison of actual 3-rank data with predicted value (\ref{Malleprobfmla}) and (\ref{AchterFmla}) over  $\finfldq{7}$}
\centering
\includegraphics[width=350pt, height=233pt]{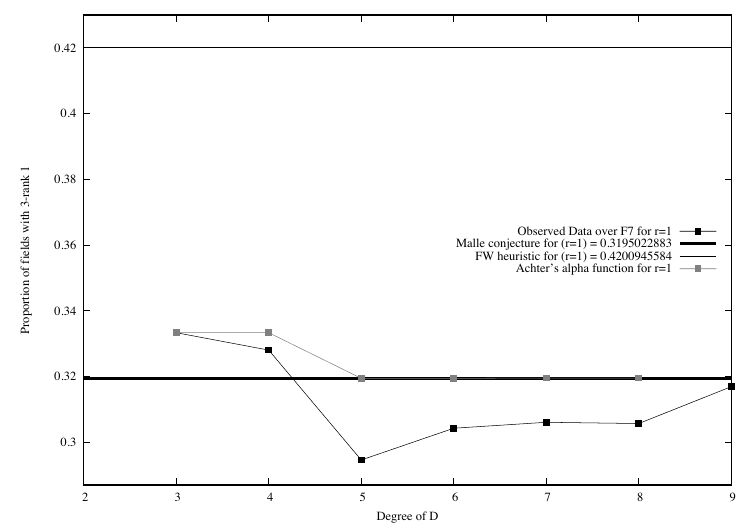}
\label{F7imagcompgraph}
\end{figure}

We also compared our $\finfldq{7}$ data to the the expression given in (\ref{FWprobfmla}).  As expected, the data is a poor fit to this value.    As noted in Malle \cite{Malle}, the Cohen-Lenstra-Martinet heuristics for $\ell$-ranks may fail when primitive $\ell$-th roots of unity lie in the base field.  This phenomenon occurs for function fields as well, in particular for the case of $\finfldq{7}$.  For $\ell=3$, we compared our data for \finfldq{7}\ to the Garton/Malle formula (\ref{Malleprobfmla}) and Achter's $\alpha(g,r)$ function (\ref{AchterFmla}) in Figure~\ref{F7imagcompgraph}.   
As seen in Figure~\ref{F7imagcompgraph}, our data agrees more closely with these predictions. The values of $\alpha(g,1)$ for $g\leq3$ appear in \cite{Achter1}.    Our data agrees more closely, and in one case exactly, with this formula.

We conclude by noting that our data compares favourably to both the Ellenberg et al.\ and Garton distributions, again giving evidence that the dependence on $q \rightarrow \infty$ could perhaps eventually be removed.

\section{Conclusion}\label{sec:concl}
This paper presented a method for finding quadratic function fields of high $3$-rank of discriminant $D$ with $-3D$ imaginary or unusual.  Our computations were carried out up to degree bounds $\deg(D) \leq 13$ for $q=5$, and various smaller bounds for $q=7,11,13$.  In addition, we very able to give an improved algorithm in the imaginary case, which allowed us to extend our computations further than those in \cite{Pieter2}.  Finally, we accumulated enough data to allow for a meaningful comparison of our $3$-rank data with the Friedman-Washington heuristics \cite{FriedWash}, along with the results of Ellenberg et al.\ \cite{EllenVenkWest}, Malle \cite{Malle,Malle2} and Garton \cite{Garton2,Garton}.  Our results approach these heuristics as $\deg(D)$ grows.

It would be interesting to modify the algorithm to the case where $-3D$ is real.  As noted in \cite{Pieter2,Pieter3}, this is currently an open problem.  In this case it is unclear how to single out efficiently a unique reduced representative in each equivalence class of binary cubic forms,
 and in fact there are generally exponentially many reduced forms in a given equivalence class.  

We were unable to make use of $\ratsfuncq{q}$-automorphisms for the case where $-3D$ is unusual.  This is mainly due to the fact that in this case, translation of a reduced form by $\beta$ results in a partially reduced form and not a reduced form.  Losing this unique form in the first part of the modified algorithm means a different approach must be used in order to maintain accurate discriminant counts.  One could attempt to just use partially reduced forms instead, but a method to avoid double-counting is needed.  Storing all the forms found so far in memory would be one solution, but would result in a significantly slower algorithm with more overhead than the original algorithm.  Trying to find specific symbolic transformation matrices to get from a type 1 form to a reduced form also seems to be a challenge.  In fact, it is unknown to the authors how many type 1, 2 or 3 forms lie in a given cycle of partially reduced unusual forms.  Getting this or any speed-up to work in the unusual case is currently an open problem.

Computing $\ell$-ranks for quadratic function fields where $\ell \neq 3$ via a similar indirect technique would also be of interest.  Unfortunately, a special connection between certain function fields and $\ell$-rank values of lower degree function fields, as given by Hasse in the case of cubic fields and $\ell=3$
has not been explored for higher degree function fields or other values of $\ell$.  Other techniques for computing $\ell$-ranks of quadratic function fields are currently being investigated.  Some recent work has been done on generating data and developing a heuristic for the case where the prime $p$ for which the $p$-rank is under consideration is equal to the characteristic of the field.  See \cite{CaisBrownEllen} for some new results along these lines.

\bibliographystyle{amsplain}
\bibliography{reffile}

\end{document}

%% file: custommacros.tex
\newtheorem{theorem}[algorithm]{Theorem}{\bfseries}{\itshape}
\newtheorem{Prop}[algorithm]{Proposition}{\bfseries}{\itshape}
\newtheorem{Lem}[algorithm]{Lemma}{\bfseries}{\itshape}
\newtheorem{Cor}[algorithm]{Corollary}{\bfseries}{\itshape}
{\bfseries}{\itshape}
\newtheorem{conj}[algorithm]{Conjecture}{\bfseries}{\itshape}
\newtheorem{definition}[algorithm]{Definition}{\bfseries}{\itshape}

\newcommand{\newU}{\ensuremath{\mathcal{U}}}
\newcommand{\sgn}{\ensuremath{\mathrm{sgn}}}
\newcommand{\finfldq}[1]{\ensuremath{\mathbb{F}_{#1}}}
\newcommand{\polyringq}[1]{\ensuremath{\mathbb{F}_{#1}[t]}}
\newcommand{\genlinfunc}{\ensuremath{GL_2(\polyringq{q})}}
\newcommand{\ints}{\ensuremath{\mathbb{Z}}}
\newcommand{\rats}{\ensuremath{\mathbb{Q}}}
\newcommand{\ratsfuncq}[1]{\ensuremath{\mathbb{F}_{#1}(t)}}

\newcommand{\ratsfuncqextn}[2]{\ensuremath{\mathbb{F}_{#1}(t,#2)}}
\newcommand{\charfld}{\ensuremath{\mathrm{char}}}